\documentclass[a4paper,11pt]{article}
\usepackage{latexsym,amssymb,amsmath}
\usepackage[all]{xy}
\CompileMatrices
\newcommand{\Z}{\mathbb{Z}}
\newcommand{\N}{\mathbb{N}}
\newcommand{\Cst}{\mathbb{C}^*}
\newcommand{\C}{\mathbb{C}}
\newcommand{\Q}{\mathbb{Q}}
\newcommand{\R}{\mathbb{R}}

\newcommand{\liea}{{\frak g}}
\newcommand{\csa}{{\frak h}}

\newcommand{\alc}{{\frak a}}
\newcommand{\loopgr}{{\cal L}G }
\newcommand{\goloopgr}{{\cal L}_{Z,\Gamma}\tilde{G}}
\newcommand{\oloopgr}{{\cal L}_{Z}\tilde{G}}
\newcommand{\looplie}{{\cal L}\liea }
\newcommand{\kmliea}{\widehat{{\cal L}\liea }^k }
\newcommand{\kmG}{\widehat{{\cal L}G }^k }

\newtheorem{theorem}{Theorem}[section]
\newtheorem{corollary}[theorem]{Corollary}
\newtheorem{lemma}[theorem]{Lemma}
\newtheorem{definition}{Definition}[section]
\newtheorem{proposition}[theorem]{Proposition}
\def\qed{\blacksquare}
\def\W{{\cal W}}

\title{A Steinberg Cross-Section for Non-Connected Affine Kac-Moody Groups}
\author{Stephan Mohrdieck\\Mathematisches Institut \\Universit\"at
  Basel\\ Rheinsprung 21\\ CH 4051 Basel\\E-mail: mohrdis@math.unibas.ch}
\begin{document}
\maketitle
%%%
%%%
%%%       Abstract
%%%
%%%
\begin{abstract}
In this paper we generalise the concept of a Steinberg
cross-section to non-connected affine Kac-Moody groups. As in the connected case, 
which was treated by G. Br\"uchert \cite{bruechert}, see also 
\cite{steinberg:regularelements,ich:conjugacyclasses} for simple algebraic groups,
a quotient map exists only on a certain 
submonoid of the affine Kac-Moody group. Non-connected affine Kac-Moody groups appear 
naturally as semidirect product of $\Cst$ with a central extension of loop 
groups ${\cal L}G$, where the underlying simple group $G$ is no longer simply 
connected and might even be non-connected, but having a simple identity component, as in \cite{ich:conjugacyclasses}. 
In contrast to the connected case, the understanding of central extensions
of non-connected loop groups  is a rather complicated issue. 
The component group $\pi_0({\cal L}G)$ is a semidirect product $\pi_1(G)\rtimes\pi_0(G)$. 
Following the approach of V. Toledano Laredo \cite{toledanolaredo:representations}, see also 
\cite{presslaysegal:loopgroups}, who dealt with the case of 
automorphisms coming from the fundamental group $\pi_1(G)$, we classify 
all of these central extensions for ${\cal L}G$ having a cyclic component group. 
Furthermore, using characters we define a quotient map for the conjugacy action on a certain submonoid of the Kac-Moody group. 
For this quotient map
we construct the cross-section in every connected component
of the Kac-Moody group  and show that, due
the $1$-dimensional centre, it carries a 
natural $\Cst$-action which does not have an analogue in the finite dimensional case. 
\end{abstract}
%%%
%%%
%%% Key words
%%%
%%%
\section{Introduction}
%%%
%%%
%%%  Motivation.
%%%
%%%
The aim of this paper is the construction of a Steinberg cross section, a section to the adjoint quotient map for
non-connected affine Kac-Moody groups. This work generalises previous results known for connected semisimple
algebraic group which were proven by Steinberg \cite{steinberg:regularelements}, see \cite{ich:conjugacyclasses}
for the non-connected case, \cite{bruechert} for connected affine Kac-Moody groups and also \cite{mokler:steinbergsection}
for a weaker result in the indefinite case.
The adjoint quotient is the quotient with respect to the conjugacy action. 
Application of the existence of the Steinberg cross section are found in the theory of principal bundles over an elliptic curve 
as well as in singularity theory, \cite{helmkeslodowy}.

%%%
%%%
%%% Description of the setting
%%%
%%%
An affine Kac-Moody group is a semidirect product $\widehat{{\cal L}G}:=\widetilde{{\cal L}G}\rtimes\Cst$ of $\Cst $ with the 
centrally extended holomorphic loop group corresponding to a simple algebraic group $G$. For simply connected $G$ 
the loop group is known to be connected and even simply connected and its central extensions are classified by its level which is an integer.
If $G$ is no longer simply connected or not even connected the situation becomes more complicated, 
see \cite{toledanolaredo:representations} where the connected but no longer simply connected case is treated. 
However, as we will show in Section 3.2, ${\cal L}G$ is a semidirect product of its identity component and its component group.
If we restrict ourselves to a cyclic subgroup $\Sigma $ of the component group and consider only central extensions of this 
subgroup then its central extension are classified by a multiple of a certain fundamental level $k_f$ (Theorem \ref{classoscentralexts}).
This restriction is justified because we are interested in the restriction of the conjugacy action of the identity component to a given exterior 
connected component.
In order to endow this central extension with a $\Cst$-action we have to consider a finite cover $\widetilde{\Cst}\to\Cst$.
To be precise we are not dealing with the loop group but with group of `open loops' as in 
\cite{toledanolaredo:representations}, ie. the group of paths from $\C$ to the universal cover
$\tilde{G}$ of $G$ whose endpoints differ by an element of the group of covering transformations. The reason is that this group 
has the right representation theory.

An induction procedure \cite{toledanolaredo:representations} provides us with the classification of all highest weight modules of 
the Kac-Moody group $\widehat{{\cal L}G}$. For the definition of the adjoint quotient map we need the characters 
$\chi_{\Lambda_0},..., \chi_{\Lambda_{s-1}},\chi_{\delta}$ of the generators $\Lambda_0,...,\Lambda_{s-1},\delta$ of 
$\Sigma $-fixed point weight lattice. As proven in \cite{etingoffrenkelkirilov,wendt:twistedcharacters} these characters 
can be defined as holomorphic and conjugacy invariant functions on the submonoid 
$\widehat{{\cal L}G}_{<1}:=\widetilde{{\cal L}G}\rtimes D^{\ast}\subset \widehat{{\cal L}G}$,
$D^{\ast}$ being the punctured unit disk.
The characters satisfy a Kac-Weyl type character formula \cite{wendt:twistedcharacters,fuchs} and
correspond to characters of a different Kac-Moody group whose Dynkin diagram is obtained from the original one by a 
folding process.

Let us mention as an aside that these characters 
and appear in a conjecture on Verlinde formulas for non-simply connected Lie groups \cite{fuchsschweigert}.

Then the adjoint quotient map $\chi$ is defined on the intersection 
$\widehat{{\cal L}G}_{<1, \tau}:= \widehat{{\cal L}G}_{<1}\cap\widehat{{\cal L}G}_{\tau}$ where $\widehat{{\cal L}G}_{\tau}$ is the 
connected component generating the component group: 
\begin{eqnarray*}
\chi:\widehat{{\cal L}G}_{<1, \tau}&\to& \C^s\times D^{\ast}\\
g&\mapsto&(\chi_{\Lambda_0}(g),..., \chi_{\Lambda_{s-1}}(g),\chi_{\delta}(g)).
\end{eqnarray*}
In general, for Kac-Moody groups the existence of an adjoint quotient, even in the affine setting, is a complicated issue; 
but, recently, there has been some progress in this direction, see \cite{mokler:adjointquotient}. 
However, dealing with holomorphic loops our notion of quotient can be justified via the moduli space of bundles over an elliptic curve.
Here we consider $G$ to be connected yielding $\Sigma=\pi_1(G)$.
For fixed $\tilde{q}\in D^{\ast}$ the space $E_{\tilde{q}}:=\Cst/\tilde{q}^{\Z}$ is an elliptic curve. The corresponding moduli space
of semistable $G$-bundles over $E_{\tilde{q}}$ has $|\Sigma|$-many connected components labelling the topological type of the bundles.
It was observed by Looijenga \cite{looijenga} that isomorphism classes of these bundles of topological type $\tau $ 
are in bijection to ${\cal L}G_0$-orbits on ${\cal L}G_{\tau}\times\{\tilde{q}\}$.
Following \cite{looijenga} and also \cite{freedmanmorganii}, respectively \cite{scweigert:moduli} for the non simply connected case,
the corresponding component of the moduli space  
is given as weighted projective space $\C^s\backslash\{0\}/\Cst$ where $\Cst $ acts like the centre of the Kac-Moody group on $\C^s$. 

We construct a local section $S:\C^s\times D^{\ast}\to \widehat{{\cal L}G}_{<1, \tau}$ to the adjoint quotient map whose image consists 
of a `twisted Coxeter cell' 
$U_{\tiny\mbox{cox}^{\tau}}\mbox{cox}^{\tau}\times D^{\ast}\subset \widehat{{\cal L}G}_{<1,\tau}$. Let us mention that for finite Coxeter groups with
outer symmetries the  twisted Coxeter element was introduced by Springer \cite{springer:regularelements}.
Our main result (Theorem \ref{maintheorem}) states that the restriction $S|_{\C^s\times\{\tilde{q}\}}$ is a section to
$\chi|_{\widetilde{{\cal L}G}\times\{\tilde{q}\}}$ for $|\tilde{q}|$ small. One of the key results for proving the theorem is the fact that 
the image of $S|_{\C^s\times\{\tilde{q}\}}$ can be endowed with a $\Cst$-action (Lemma \ref{Cstarlemma}) which has no analogue in the finite 
dimensional situation, however see \cite{bruechert} for connected affine Kac-Moody groups.

Using this section and its $\Cst$-action gives a new proof of the result mentioned above that the moduli space 
of semistable bundles over $E_{\tilde{q}}$ in each component is a weighted projective space. These ideas will be published elsewhere, \cite{wendt:bundles}. 

Furthermore, the understanding of conjugacy classes and the quotient map is used by Slodowy and Helmke \cite{helmkeslodowy} to 
establish a connection between simple elliptic singularities and Kac-Moody groups.

%%%
%%%
%%%  Organisation of paper
%%%
%%%
The paper is organised as follows: In Section 2 we fix our notation and review the folding process of Dynkin diagrams.
Furthermore, we introduce the notion of twisted Coxeter elements and investigate its properties on the standard representation of the 
corresponding Coxeter group. The affine Kac-Moody groups and Lie algebras are introduced in Section 3. Here, we will discuss in detail the classification
of central extension for non-connected loop groups and summarise their representation theory.
The definition of the adjoint quotient map and construction recipe for the section will be contained in Section 4.
Furthermore, in this section we will formulate and prove the main result.
In the appendix we compile table containing several data appearing in this paper.

%%%
%%%
%%% Acknowledgement
%%%
%%%
{\it Acknowledgements:} The author would like to thank the Schweizerischer Nationalfonds for financial support.
%%%
%%%
%%%  Section 2: Symmetries of Affine Dynkin Diagrams
%%%
%%%
\section{Symmetries of Affine Dynkin Diagrams}
%%%
%%%
%%%  Subsection 2.1: Basic definitions
%%%
%%%
\subsection{Basic definitions}
Let $\Pi $ be a Dynkin diagram of finite type of rank $r$ with Cartan matrix $C$. Furthermore, denote by $\Delta $ the corresponding 
root system, by $\check{\Delta}$ its dual, and by ${\cal W}$ its Weyl group. We use the symbol $\Pi:=\{\alpha_1,...,\alpha_r\}$ 
also for the set of simple roots. The corresponding simple reflections are denoted by $s_i,i\in\{1,..,r\}$.
The integral span $Q:=\Z\Delta $ is the called the root lattice and $\check{Q}:=\Z\check{\Delta}$ the co-root lattice.
The corresponding dual lattices $\check{P}:=\mbox{Hom}(Q,\Z)$ respectively $P:=\mbox{Hom}(\check{Q},\Z)$ are the 
co-weight lattice respectively the weight lattice. The highest element of $\Delta $ will be denoted by $\theta $.

Our next step is to describe briefly the affine untwisted root systems. For details of the construction we refer to Kac's book
\cite{kac:infdimliealgebras}. 
The affine untwisted diagram is obtained by adding an additional vertex and some lines to the corresponding 
finite dimensional diagram. We write $\widehat{\Pi} $ for the corresponding Dynkin diagram as well as for the set of simple roots 
$\{\alpha_0,...,\alpha_r\}$ and $\widehat{C}$ for its Cartan matrix. 
The set of simple roots generates a free complex vector space $\csa^{'*}$. 
Furthermore we introduce the set of simple co-roots $\check{\widehat{\Pi}}:=\{\check{\alpha}_0,..,\check{\alpha}_r\}$
spanning its dual $\csa'$ by the rule $ \check{\alpha}_i(\alpha_j)=\widehat{C}_{ij}$
Note that the Cartan matrix has a one-dimensional kernel $\Z\delta $ with $\delta:=a_0\alpha_0+...+a_r\alpha_r$ 
where the $a_i$ can be chosen to be positive integers, the {\it Kac-labels}. 
Then there is the equality $\alpha_0=\delta-\theta $.
In terms of the underlying finite dimensional root system $\Delta$ the affine root system $\widehat{\Delta}$ has the following description:
\begin{equation}
\widehat{\Delta}:=\{\alpha+n\delta, \alpha\in\Delta \mbox{and}\, n\in\Z\}\cup\Z\backslash\{0\}\delta .
\end{equation}
It decomposes into the disjoint union of the set of {\it real roots}  $\widehat{\Delta}^{re}:=\widehat{\Delta}\backslash(\Z\backslash\{0\})\delta$ 
and the set of {\it imaginary roots} $\widehat{\Delta}^{im}:=\Z\backslash\{0\}\delta$.
In analogy to the finite dimensional situation the simple reflections $s_i,i\in\{0,..,r\}$ corresponding to the set of simple root generate the 
affine Weyl group $\widehat{\W}$. Here $s_i$ is defined by $s_i(\alpha_j)=\alpha_j-\widehat{C}_{ij}\alpha_i$ yielding a linear action of 
$\widehat{\W}$ on $\csa^{'*}$ and dually on $\csa' $.
Then, we have $\widehat{\Delta}^{re}=\widehat{\W}\widehat{\Pi}$.
In terms of the finite Weyl group there is the identity $\widehat{\W}=\check{Q}\rtimes\W $.
Here, an element $\check{\beta}\in\check{Q}$ acts on $\R\otimes\check{Q}$ as translation $\tau_{\check{\beta}}$
and $s_0 $ is identified with $s_{\theta}\circ\tau_{-\check{\theta}} $, where $\check{\theta}\in\check{\Delta}$ is the co-root corresponding to $\theta$. 
A fundamental domain for the action of this group on $\R\otimes\check{Q}$ is given by the fundamental alcove 
$\alc:=\{x\in\R\otimes\check{Q}|\alpha_i(x)\geq 0, i\in\{1,...,r\}\,\mbox{and}\,\theta(x)\leq 1\}$.
Let $\check{\chi}$ be any $\W $-stable lattice with $\check{Q}\subset\check{\chi} \subset\check{P}$. (Hence, $\check{\chi}$ is the co-character lattice
of a simple algebraic group with Dynkin diagram $\Pi$.) Then we define the group $\widehat{\W}(\check{\chi}):=\check{\chi}\rtimes\W $.
The group $\widehat{\W}(\check{P})$ is called the {\it extended affine Weyl group}.
As is shown in \cite{bourbaki:lie} the stabiliser $\widehat{\W}(\check{\chi})_{\alc}$ of $\alc$ in $\widehat{\W}(\check{\chi})$ is 
mapped isomorphically to $\check{\chi}/\check{Q}=\widehat{\W}(\check{\chi})/\widehat{\W}$.
Thus there is a semidirect product decomposition:
\begin{equation}
\widehat{\W}(\check{\chi})=\widehat{\W}\rtimes\widehat{\W}(\check{\chi})_{\alc}.
\label{decompossemidir}
\end{equation}
Since the vertices of $\alc$ correspond to the vertices of $\widehat{\Pi }$
the group $\widehat{\W}(\check{\chi})_{\alc}$ is a subgroup of the diagram automorphisms $\mbox{Aut}(\widehat{\Pi})$
of the affine Dynkin diagram.

Our next goal is to give a description of the whole group $\mbox{Aut}(\widehat{\Pi})$.
This can be defined as the group of permutations $\tau $ of the index set $\{0,...,r\}$ satisfying $\widehat{C}_{\tau(i)\tau(j)}=
\widehat{C}_{ij}$ for all $i,j\in\{0,...,r\}$.
The following lemma is known in Lie theory:
\begin{lemma} 
\label{descriptionauts}
The automorphism group $\mbox{Aut}(\widehat{\Pi})$ is a semidirect product:
\begin{equation}
\mbox{Aut}(\widehat{\Pi})=\check{P}/\check{Q}\rtimes\mbox{Aut}(\Pi)=\widehat{\W}(\check{P})_{\alc}\rtimes\mbox{Aut}(\Pi).
\end{equation}
Here $\mbox{Aut}(\Pi)$ is the group of symmetries of the Dynkin diagram of finite type $\Pi$.
\end{lemma}
As a consequence we get the statement:
\begin{corollary}
For the whole automorphism group of the affine root system we get:
\begin{equation}
\mbox{Aut}(\widehat{\Delta})=\widehat{\W}\rtimes\mbox{Aut}(\widehat{\Pi})=
\widehat{\W}\rtimes(\check{P}/\check{Q}\rtimes\mbox{Aut}(\Pi))
=\widehat{\W}(\check{P})\rtimes\mbox{Aut}(\Pi). 
\end{equation}
\end{corollary} 

Our next goal is to define the root and the weight lattice associated to affine Dynkin diagrams. 
First consider the complex vector space $\csa $ defined by $\csa :=\csa'\oplus\C d$.  
Here $d$ is an additional variable and we extend the linear forms $\alpha_i$ to $\csa $ by:
\begin{eqnarray}
\alpha_i(d)&:=&0,\quad i\in\{1,...,r\}\\
\delta(d)&:=&1.
\end{eqnarray}
The $\widehat{\W}$-action extends to $\csa $ by defining $s_i(d):=d-\alpha_i(d)\check{\alpha}_i$.
Furthermore we can lift the $\mbox{Aut}(\widehat{\Pi})$-action to $\csa $ by requiring $\alpha_{\tau(i)}(d)=\alpha_i(\tau(d))$
for any $\tau\in\mbox{Aut}(\widehat{\Pi})$.
Denote by $\csa^*$ the dual of $\csa$. 
The affine co-root lattice is defined by $\check{\widehat{Q}}:=\Z<\check{\alpha}_0,...,\check{\alpha}_r,d>$. 
Its dual $\widehat{P}:=\mbox{Hom}(\check{\widehat{Q}},\Z)$ is called the weight lattice. We write $\{\lambda_0,...,\lambda_r,\delta\}$
for the corresponding dual basis. The cone $\widehat{P}^+:=\N_0<\lambda_0,...,\lambda_r>\oplus\Z\delta $
is called the {\it cone of dominant weights}.  

As conclusion of this subsection let us briefly describe the case of twisted affine root systems. 
Let $\sigma\in\mbox{Aut}(\Pi)$ and denote by $\csa^{\sigma}$ the fixed point space $\{h\in\csa,\,\sigma (h)=h\}$.
The twisted affine Dynkin diagram $\widehat{\Pi}(\sigma)$ 
can be described as $\widehat{\Pi}(\sigma)=\Pi^{\sigma}\cup\{\delta-\theta'\}$.
Here $\Pi^{\sigma}$ is the Dynkin diagram of the root system $\{\alpha|_{\csa^{\sigma}},\,\alpha\in\Delta\}$,
except for the case $\mbox{A}_{2n}$,
and $\theta'$ is the highest short root. Correspondingly, there is the Cartan matrix $\widehat{C}(\sigma)$.
In the $\mbox{A}_{2n}$-case the root system $\{\alpha|_{\csa^{\sigma}},\,\alpha\in\Delta\}$ is of type $\mbox{BC}_n$.
Here we define $\Pi^{\sigma}$ to be a basis of the subsystem of type $\mbox{C}_n$ containing all roots which  are long or of intermediate length.
The twisted affine Weyl group is given by $\widehat{\W}(\sigma)=\W^{\sigma}\ltimes Q^{\sigma}$, where 
$\W^{\sigma}=\{w\in\W, \,\sigma w=w\sigma\}$ is the fixed point Weyl group and $Q^{\sigma}$ the lattice generated by $\Pi^{\sigma}$.
Furthermore, for any $\sigma$-stable character lattice $\chi$ there is the extended Weyl group 
$\widehat{\W}(\chi)(\sigma)=\W^{\sigma}\ltimes\chi^{\sigma}$. In this situation the notion of fundamental 
alcove $\alc$ does exist as well, and Equation \ref{decompossemidir} still holds.
It turns out that the only interesting cases are $\chi=Q,P$. 
As analogue of Lemma \ref{descriptionauts} we get:
\begin{equation}
\mbox{Aut}(\widehat{\Pi}(\sigma))=P^{\sigma}/Q^{\sigma}=\widehat{\W}(P)(\sigma)_{\alc}. 
\end{equation}
We define the weight lattice and the cone of dominant weights analogously.

In the sequel $\widehat{C}$ denotes any affine Cartan matrix, twisted or untwisted.

%%%
%%%
%%%  Subsection 2.2: The Symmetries of the affine Dynkin Diagrams
%%%
%%%
\subsection{The Symmetries of the affine Dynkin Diagrams}
Here, we give a list of the symmetry groups of all affine Dynkin diagrams, Furthermore, we choose representatives of 
their conjugacy classes which we use in the sequel. We label the vertices of the Dynkin diagram as in \cite{bruechert}.

{\bf $\mbox{A}^1_1$:} The group of diagram automorphisms is $\mbox{Aut}(\widehat{\Pi})=\check{P}/\check{Q}=\Z/2\Z$ and we denote 
by $\gamma $ its generator.

{\bf $\mbox{A}^1_n,\,n\geq 2$:} In this case we have $\mbox{Aut}(\widehat{\Pi})={\cal D}_{n+1}$, the dihedral group containing
$2n+2$ elements. We denote by $\sigma $ the generator of $\mbox{Aut}(\Pi)=\Z/2\Z$ and by $\gamma $ a generator of the cyclic part
$\check{P}/\check{Q}=\Z/(n+1)\Z$. These elements can be described explicitly by: $\sigma(0)=0$ and $\sigma (i)=n+1-i$ respectively 
$\gamma (i)=i+1\,\mbox{mod}\,n+1$.
We choose the following representatives of ${\cal D}_{n+1}$-conjugacy classes:\\
{\it $n$ even:} $\gamma^l,\,l|n+1$ and $\sigma$.\\
{\it $n$ odd:} $\gamma^l,\,l|n+1$, $\sigma$, $\sigma\gamma$.

{\bf $\mbox{B}^1_n,\,n\geq 3$:} The group of diagram automorphisms is $\mbox{Aut}(\widehat{\Pi})=\check{P}/\check{Q}=\Z/2\Z$ 
with generator $\gamma $.

{\bf $\mbox{C}^1_n,\,n\geq 2$:} Here we have $\mbox{Aut}(\widehat{\Pi})=\check{P}/\check{Q}=\Z/2\Z$ 
with generator $\gamma $.

{\bf $\mbox{D}^1_4$:} The automorphism group $\mbox{Aut}(\widehat{\Pi})$ is the symmetric group ${\cal S}_4$.
In this case we have $\check{P}/\check{Q}=\Z/2\Z\times\Z/2\Z$ and 
$\mbox{Aut}(\Pi)={\cal S}_3$. ${\cal S}_4$ is known to have four non-trivial conjugacy classes.
We choose generators $\rho $ respectively $\sigma $ of order $3$ 
respectively $2$ of ${\cal S}_3$ acting on the diagram by: $\rho(1)=3$, $\rho(3)=4$, $\rho(4)=1$ fixing the
vertices $0$ and $2$ and $\sigma $ by interchanging the vertices $3$ and $4$ fixing the other three.
Furthermore, we define the element $\gamma $ defined by $\gamma(0)=4$, $\gamma(1)=3$, $\gamma(2)=2$, 
$\gamma(3)=0$ and $\gamma(4)=1$ of order $4$. The missing conjugacy class is represented by $\gamma^2$.
Then $\check{P}/\check{Q}$ is given by $\{id,\,\gamma^2 ,\,\gamma\sigma ,\,\gamma^3\sigma \}$.

{\bf $\mbox{D}^1_n,\,n\geq 5$:} In this case we have $\mbox{Aut}(\widehat{\Pi})={\cal D}_{4}$ and
$\mbox{Aut}(\Pi)=\Z/2\Z$. Choose the following elements as representatives for the conjugacy classes:
Let $\sigma $ be the generator of $\mbox{Aut}(\Pi)$ which is defined by interchanging the vertices 
with label $n-1$ and $n$ and choose the element $\gamma$ of order $4$ given by $\gamma(0)=n$, $\gamma(1)=n-1$, 
$\gamma(n-1)=0$ and $\gamma(n)=1$, fixing all remaining  vertices. Then the missing conjugacy classes are 
represented by $\gamma^2$ and $\gamma\sigma$. For identifying the elements of $\check{P}/\check{Q}$
we have to distinguish two cases:\\
{\it $n$ odd:} Then $\check{P}/\check{Q}=\Z/4\Z=\{id,\,\gamma,\gamma^2,\,\gamma^3\}$.\\
{\it $n$ even:} Then $\check{P}/\check{Q}=\Z/2\Z\times\Z/2\Z=\{id,\,\gamma\sigma,\gamma^2,\,\gamma^3\sigma\}$.

{\bf $\mbox{E}^1_6$:} Here we have the symmetry group $\mbox{Aut}(\widehat{\Pi})={\cal S}_{3}$ with 
$\check{P}/\check{Q}=\Z/3\Z$ respectively $\mbox{Aut}(\Pi)=\Z/2\Z$. We choose as representatives of he nontrivial
conjugacy classes the generator $\sigma $ of $\mbox{Aut}(\Pi)$ and a generator $\gamma$ of $\check{P}/\check{Q}$.
These elements are explicitly given by  given by: $\sigma(1)=5$, $\sigma(2)=4$, $\sigma(4)=2$, $\sigma(5)=1$ with 
$\sigma$ fixing the remaining vertices and $\gamma(0)=1$, $\gamma(1)=5$, $\gamma(2)=4$, $\gamma(3)=3$, $\gamma(4)=6$, 
$\gamma(5)=1$ and $\gamma(6)=0$.

{\bf $\mbox{E}^1_7$:} The group of diagram automorphisms is $\mbox{Aut}(\widehat{\Pi})=\check{P}/\check{Q}=\Z/2\Z$ 
with generator $\gamma $.

{\bf $\mbox{A}^2_{2n-1},\,n\geq 3$:} Here $\mbox{Aut}(\widehat{\Pi})$ is given by $\Z/2\Z$ 
with generator $\gamma $.

{\bf $\mbox{D}^2_{n=1},\,n\geq 2$:} The group of diagram automorphisms is $\mbox{Aut}(\widehat{\Pi})=\Z/2\Z$ 
with generator $\gamma $.
%%%
%%%
%%%  Subsection 2.3: Folding of Dynkin Diagrams
%%%
%%%
\subsection{Folding of Dynkin Diagrams}
An important tool for the understanding of the representation theory of the non-connected affine Kac-Moody
groups is the folding of Dynkin diagrams. This was considered first by Jantzen \cite{jantzen:darstellungen} in the finite dimensional case and 
by Fuchs, Schweigert and Schellekens \cite{fuchs} for affine Kac-Moody algebras.
Let us summarise the construction recipe here:

To a pair $(\widehat{C},\,\tau )$ consisting of an affine (twisted or untwisted) Cartan matrix $\widehat{C}$ and a symmetry $\tau\in\mbox{Aut}(\widehat{\Pi})$ thereof 
we associate another Cartan matrix ${}^{\tau }\widehat{C}$ according to the following rule:
Denote by $\Sigma $ the cyclic group of $\mbox{Aut}(\widehat{\Pi})$ generated by $\tau $ and assume $\tau $ to be of order $N$.
Using $I$ for the index set $\{0,...,r\}$ set $t_i:=|\{j\in I|\,\widehat{C}_{ij}\neq 0\}| $, the number of all elements in the $\Sigma $-orbit of $i$ which
are adjacent to $i$ (including $i$ itself). Taking a look at the Cartan matrices of affine type we see:
$t_i\leq 2$ except for the case $\mbox{A}^1_n$ with $\tau $ being the cyclic permutation of order $n+1$ where
we have $t_i=3$.
Now, we define the matrix  $\widetilde{C}$ by taking for its columns the sum over all $\Sigma $-orbits of the columns of $\widehat{C}$
and multiplying with $t_i|\Sigma _i|$: (Here $\Sigma _i$ is the stabiliser of $i\in I$ in $\Sigma $.)
%%%
%%%
%%% Definition 2.1: Folded Cartan matrix and orbit Lie algebra
%%%
%%%
\begin{equation}
\label{eqnfoldeddiagram}
\widetilde{C}_{ij}:=t_j|\Sigma_j|\sum_{k=0}^{N-1}\widehat{C}_{i\tau^k(j)}.
\end{equation}
\begin{definition}
The matrix obtained from $\widetilde{C}$ by removing redundant columns and rows is called the {\it folded Cartan matrix}
and denoted by ${}^{\tau }\widehat{C}$. It is of square size and has as many rows as there are $\Sigma $-orbits on the 
index set $I$. Its Dynkin diagram is called the folded Dynkin diagram ${}^{\tau }\widehat{\Pi}$.
\end{definition}

The folded Dynkin diagrams associated to all pairs $\widehat{C},\tau $ can be found in the appendix, Table 1.

{\it Remark:} We can give the following pictorial description of the folding process:
Take a Dynkin diagram $\widehat{\Pi }$ together with the diagram automorphism $\tau $.
As above, let $\Sigma $ be the subgroup of $\mbox{Aut}(\widehat{\Pi })$ generated by $\tau $.
The vertices of the folded diagram correspond to the $\Sigma $-orbits of vertices of $\widehat{\Pi} $.
Let $i$ and $j$ be two distinct $\Sigma$-orbits and denote by $|\Sigma_i|$ respectively $|\Sigma_j|$ the order
of each of its elements. We can  assume $|\Sigma_i|\geq |\Sigma_j|$. In general, connect the vertices $i$ and $j$ 
by $\frac{|\Sigma_i|}{|\Sigma_j|}\mbox{max}_{p\in i,\,q\in j}|\{\mbox{lines between $p$ and $q$ in $\widehat{\Pi}$}\}|$-many lines
keeping any existing arrows. In case there do not exist any arrows between the representatives in $\widehat{\Pi}$ we add one 
pointing from $j$ to $i$.
The only exception to this rule is the following case: $|\Sigma_i|=|\Sigma_j|$, and there is a vertex $p\in i$ 
adjacent to another vertex in the same orbit but no element $q\in j$ is adjacent to any other element in $j$.
Then we have to connect $i$ and $j$ with $2\,\mbox{max}_{p\in i,\,q\in j}|\{\mbox{lines between $p$ and $q$ in $\widehat{\Pi}$}\}|$-many
lines, keep any existing arrows and add another arrow pointing from $i$ to $j$.  

Let us illustrate this procedure by the example $\mbox{C}^1_{n}$:
%%%
%%%
%%%  Example of folding: C^1_n case
%%%
%%%
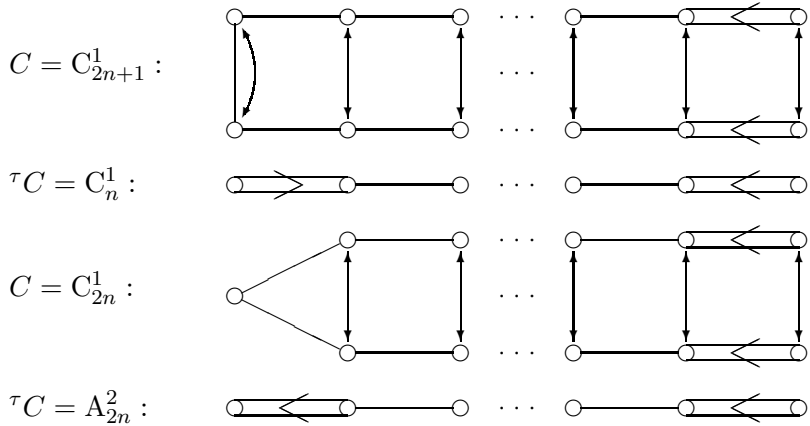
\begin{center}
%%
%% n odd
%%
\setlength{\unitlength}{1mm}
\begin{picture}(120,22)
\thinlines 
\put(0,8.5){$ C=\mbox{C}^1_{2n+1}:$} 
\put(30,1){\circle{2}} \put(60,1){\circle{2}} \put(45,1){\circle{2}} \put(75,1){\circle{2}} \put(90,1){\circle{2}}\put(105,1){\circle{2}}
\put(31,1){\line(1,0){13}} \put(46,1){\line(31,0){13}}\put(30,2){\line(0,1){13}} \put(76,1){\line(1,0){13}} \put(90,0){\line(1,0){15}}
\put(90,2){\line(1,0){15}}
\put(96,1){\line(2,1){3.75}}
\put(96,1){\line(2,-1){3.75}}
\put(45,5){\line(0,1){7}} \put(75,5){\line(0,1){7}} \put(90,5){\line(0,1){7}}\put(105,5){\line(0,1){7}}\put(60,5){\line(0,1){7}}
\put(45,5.5){\vector(0,-1){3}} \put(45,11.5){\vector(0,1){3}} \put(75,11.5){\vector(0,1){3}}\put(75,5.5){\vector(0,-1){3}} \put(90,11.5){\vector(0,1){3}}
\put(90,5.5){\vector(0,-1){3}} \put(60,5.5){\vector(0,-1){3}} \put(60,11.5){\vector(0,1){3}}\put(105,5.5){\vector(0,-1){3}} \put(105,11.5){\vector(0,1){3}}
\put(65.5,1){\circle*{.05}} \put(67.5,1){\circle*{.05}} \put(69.5,1){\circle*{.05}}
\put(30,16){\circle{2}} \put(60,16){\circle{2}} \put(45,16){\circle{2}} \put(75,16){\circle{2}} \put(90,16){\circle{2}}\put(105,16){\circle{2}}
\put(31,16){\line(1,0){13}}\put(46,16){\line(1,0){13}} \put(76,16){\line(1,0){13}} \put(90,15){\line(1,0){15}}\put(90,17){\line(1,0){15}}
\put(96,16){\line(2,1){3.75}}
\put(96,16){\line(2,-1){3.75}}
\qbezier(31.2,3)(34,8.5)(31.2,14)
\put(31,2.5){\vector(-1,-2){0}}
\put(31,14.5){\vector(-1,2){0}}
\put(65.5,16){\circle*{.05}}\put(67.5,16){\circle*{.05}}\put(69.5,16){\circle*{.05}}
\put(65.5,8.5){\circle*{.05}}\put(67.5,8.5){\circle*{.05}}\put(69.5,8.5){\circle*{.05}}
\end{picture} 
%%
%% folded type: C^1_n
%%
\setlength{\unitlength}{1mm}
\begin{picture}(120,7)
\thinlines\put(0,0){${}^{\tau }C=\mbox{C}^1_n:$}
\put(30,1){\circle{2}} \put(60,1){\circle{2}} \put(45,1){\circle{2}} \put(75,1){\circle{2}} \put(90,1){\circle{2}} \put(105,1){\circle{2}}
\put(30,0){\line(1,0){15}}\put(30,2){\line(1,0){15}} \put(46,1){\line(1,0){13}}\put(76,1){\line(1,0){13}} \put(90,0){\line(1,0){15}}\put(90,2){\line(1,0){15}}
\put(39,1){\line(-2,1){3.75}}
\put(39,1){\line(-2,-1){3.75}}
\put(96,1){\line(2,1){3.75}}
\put(96,1){\line(2,-1){3.75}}
\put(65.5,1){\circle*{.05}}\put(67.5,1){\circle*{.05}}\put(69.5,1){\circle*{.05}}
\end{picture}
%%
%% n even 
%%
\setlength{\unitlength}{1mm}
\begin{picture}(120,22)
\thinlines 
\put(0,8.5){$ C=\mbox{C}^1_{2n}:$} 
\put(30,8.5){\circle{2}} \put(60,1){\circle{2}} \put(45,1){\circle{2}} \put(75,1){\circle{2}} \put(90,1){\circle{2}}\put(105,1){\circle{2}}
 \put(46,1){\line(1,0){13}} \put(76,1){\line(1,0){13}} \put(90,0){\line(1,0){15}}
\put(90,2){\line(1,0){15}}
\put(96,1){\line(2,1){3.75}}
\put(96,1){\line(2,-1){3.75}}
\put(45,5){\line(0,1){7}} \put(75,5){\line(0,1){7}} \put(90,5){\line(0,1){7}}\put(105,5){\line(0,1){7}}\put(60,5){\line(0,1){7}}
\put(45,5.5){\vector(0,-1){3}} \put(45,11.5){\vector(0,1){3}} \put(75,11.5){\vector(0,1){3}}\put(75,5.5){\vector(0,-1){3}} \put(90,11.5){\vector(0,1){3}}
\put(90,5.5){\vector(0,-1){3}} \put(60,5.5){\vector(0,-1){3}} \put(60,11.5){\vector(0,1){3}}\put(105,5.5){\vector(0,-1){3}} \put(105,11.5){\vector(0,1){3}}
\put(65.5,1){\circle*{.05}} \put(67.5,1){\circle*{.05}} \put(69.5,1){\circle*{.05}}
 \put(60,16){\circle{2}} \put(45,16){\circle{2}} \put(75,16){\circle{2}} \put(90,16){\circle{2}}\put(105,16){\circle{2}}
\put(46,16){\line(1,0){13}} \put(76,16){\line(1,0){13}} \put(90,15){\line(1,0){15}}\put(90,17){\line(1,0){15}}
\put(96,16){\line(2,1){3.75}}
\put(96,16){\line(2,-1){3.75}}
\put(30.9,9){\line(2,1){13}}\put(30.9,8){\line(2,-1){13}}
\put(65.5,16){\circle*{.05}}\put(67.5,16){\circle*{.05}}\put(69.5,16){\circle*{.05}}
\put(65.5,8.5){\circle*{.05}}\put(67.5,8.5){\circle*{.05}}\put(69.5,8.5){\circle*{.05}}
\end{picture} 
%%
%%  folded type: A^2_2n
%%
\setlength{\unitlength}{1mm}
\begin{picture}(120,7)
\thinlines\put(0,0){${}^{\tau }C=\mbox{A}^2_{2n}:$}
\put(30,1){\circle{2}} \put(60,1){\circle{2}} \put(45,1){\circle{2}} \put(75,1){\circle{2}} \put(90,1){\circle{2}} \put(105,1){\circle{2}}
\put(30,0){\line(1,0){15}}\put(30,2){\line(1,0){15}} \put(46,1){\line(1,0){13}}\put(76,1){\line(1,0){13}} \put(90,0){\line(1,0){15}}\put(90,2){\line(1,0){15}}
\put(36,1){\line(2,1){3.75}}
\put(36,1){\line(2,-1){3.75}}
\put(96,1){\line(2,1){3.75}}
\put(96,1){\line(2,-1){3.75}}
\put(65.5,1){\circle*{.05}}\put(67.5,1){\circle*{.05}}\put(69.5,1){\circle*{.05}}
\end{picture}
\end{center}
%%%
%%%
%%%  Subsection 2.4: Twisted Coxeter Elements
%%%
%%%
\subsection{Twisted Coxeter Elements}
In this subsection we define the notion of a twisted Coxeter element which was introduced by Springer \cite{springer:regularelements} 
for finite Coxeter groups
and discuss some of its properties.
The twisted Coxeter element plays a key role in our investigation.

As above let $\widehat{\Pi}$ be an affine Dynkin diagram with Cartan matrix $\widehat{C}$.
Let $I:=\{0,...,r\}$ be the index set labelling the vertices of $\widehat{\Pi}$.

Fix any element $\tau\in\mbox{Aut}(\widehat{\Pi})$ acting with $s$ orbits on the Dynkin diagram. 
After relabelling the index set we can assume that the set $\{0,..,s-1\}$ is a set of representatives of $\tau$-orbits on $I$. 

\begin{definition} 
The element $\mbox{cox}^{\tau}:=s_0...s_{s-1}\tau\in \widehat{\W}\rtimes \mbox{Aut}(\widehat{\Pi})$ is called a twisted Coxeter 
element of $\widehat{\W}\rtimes \mbox{Aut}(\Pi)$ corresponding to $\tau$.
\end{definition}
%%%
%%%
%%%  Remarks
%%%
%%%
{\it Remark:} $(i)$ Observe that we get back the usual definition of a Coxeter element for $\tau=e$.

$(ii)$ Twisted Coxeter elements can be defined for arbitrary, not necessarily affine, generalised Cartan matrices.
Lemma \ref{coxunique}, Proposition \ref{actioncartan} and Corollary \ref{fixpointcartan} still hold in this more general context.

Applying Lemma 7.5 of \cite{springer:regularelements}, see also Lemma 1, p.117 of \cite{bourbaki:lie} to our situation yields:
%%%
%%%
%%%  Lemma 2.1: Uniqueness of twisted Coxeter element
%%%
%%%
\begin{lemma}
\label{coxunique}
If the Dynkin diagram $\widehat{\Pi} $ contains no cycles, the twisted Coxeter element $\mbox{cox}^{\tau }$
is unique up to conjugation with elements in $\widehat{\W }$.
\end{lemma}
{\it Remark:} In case of the Dynkin diagram $\mbox{A}_r^1$ and $\tau $ the cyclic permutation of order $r+1$ 
it is easily verified that the corresponding twisted Coxeter element is also unique up to conjugation in 
$\widehat{\W }$, eg. $s_1\tau=s_1(s_0\tau)s_1$.
But for $\sigma =\tau^l $ with $l|r+1 $, we can indicate a counter-example:

Consider the case $\mbox{A}^1_5$. For the square $\tau^2 $ of the cyclic permutation $\tau\in\mbox{Aut}(\widehat{\Pi})$, 
$\tau(i)=i+1\,\mbox{mod}5$ we have for example the following two twisted Coxeter elements:
$(s_0\tau )^2=s_0s_1\tau^2 $ and $s_0s_3\tau^2$. The characteristic polynomials for their action on 
$\csa' $ are $(t-1)(t^5-1)$ resp. $(t^2-1)(t^4-1)$. Hence, they cannot be conjugate under $\widehat{\W} $.

In contrast to the finite dimensional case (twisted) Coxeter elements are no longer of finite order in the infinite
dimensional case. They turn even out to be no longer diagonalisable as we will show below. (See eg. Br\"uchert in case of usual Coxeter
elements in the affine case.)

As above, denote by $\Sigma $ the cyclic subgroup generated by $\tau $ in $\mbox{Aut}(\widehat{\Pi})$.
Our next aim is to investigate the action of $\mbox{cox}^{\tau }$ on the vector spaces $\csa $ and $\csa'$.
In particular we are interested in the multiplicity of eigenvalue $1$. We have a first result,
which was proven for $\tau =id $ by Coleman \cite{coleman:coxetertrafos}, Theorem 3.1: 
%%%
%%%
%%%  Proposition 2.1: Description of eigenvalue 1 spaces of certain group elements
%%%
%%%
\begin{proposition}
\label{actioncartan}
Let $\tau\in \mbox{Aut}(\widehat{\Pi})$ be a diagram automorphism and $s_{i_1},...,s_{i_p}$ simple reflections
with the $i_j$ lying in distinct $\Sigma $-orbits on $I$.  
Then the following is true for any $x\in\csa $ or $x\in\csa'$:
\begin{eqnarray}
s_{i_1}...s_{i_p}\tau(x)=x \iff & s_{i_j}(x)=x \\
&\tau(x)=x. 
\end{eqnarray}  
\end{proposition}
%%%
%%%  Proof of Proposition 2.1.
%%%
{\it Proof:} Clearly, the right hand side implies the left one. For the other direction, note that
 we have the following equivalence:
\begin{eqnarray}
s_{i_1}...s_{i_p}\tau(x)&=&x\qquad\iff\\
\tau(x)-x&=& s_{i_p}...s_{i_1}(x)-x.
\end{eqnarray}  
Now, the right hand side can be written as: 
\begin{equation}
s_{i_p}...s_{i_1}(x)-x=\sum_{j=1}^{p}c_j\check{\alpha}_{i_j}.
\end{equation}
Thus we obtain: 
\begin{equation}
\tau(x)=x+\sum_{j=1}^{p}c_j\check{\alpha}_{i_j}.
\end{equation}
Taking the sum over the whole $\Sigma $-orbit of this equation yields:
\begin{equation}
\sum_{j=1}^{p}d_j\sum (\Sigma \mbox{-orbit of }\check{\alpha}_{i_j})=0,
\end{equation} 
for $d_j=c_j |\Sigma_{\check{\alpha}_{i_j}}|$. 
Now, the condition on $i_1,...,i_p $ yields $c_j =0 $. 
This implies:
\begin{equation}
\C\check{\alpha}_{i_p}+x\ni s_{i_p}(x)=s_{i_{p-1}}...s_{i_1}(x)\in x+\bigoplus_{j=1}^{p-1}\C \check{\alpha_{i_j}}.
\end{equation}
Since the affine spaces on the right and left intersect only in $x$, 
we get $s_{i_p}(x)=x $. Now, the statement follows by induction.$\qed$ .
%%%
%%%
%%%  Corollary 2.1: Description of fixed points of twisted Coxeter element
%%%
%%%
\begin{corollary}
\label{fixpointcartan}
Let $\tau $ be as above. Denote by $\csa^{\tau }$ the fixed point subspace of $\csa $ under $\tau $.
Then  the following result for the twisted Coxeter element $\mbox{cox}^{\tau }$ holds:
\begin{equation}
\{x\in\csa , \mbox{cox}^{\tau }(x)=x\}=\csa^{\tau }\cap\bigcap_{i=0}^{s}(\mbox{ker}\,\alpha_i|_{\csa^{\tau }}).
\end{equation}
A similar result is true for $\csa '$.
\end{corollary}
Since $\mbox{dim}\csa=\mbox{dim}\csa'+1=r+2$ we conclude: 
%%%
%%%
%%%  Corollary 2.2: Restrictions of roots to fixed point subspaces
%%%
%%%
\begin{lemma}
i) $\mbox{dim}\,\csa^{\tau}=\mbox{dim}\,\csa ^{'\tau }+1=s+1$.\\
ii) The restrictions $\alpha_0|_{\csa ^{\tau }},...,\alpha_{s-1}|_{\csa^{\tau }}$ are linearly independent while the 
restrictions $\alpha_0|_{\csa ^{'\tau }},...,\alpha_{s-1}|_{\csa ^{'\tau }}$ satisfy the one equation:
\begin{equation}
a_0|\Sigma_{\alpha_0}|\alpha_0|_{\csa ^{'\tau }}+...+a_{s-1}|\Sigma_{\alpha_{s-1}}|\alpha_{s-1}|_{\csa ^{'\tau }}=0.
\end{equation}
Here the $a_i$ are the Kac labels.
\end{lemma}  
%%%
%%%
%%%  Corollary 2.3: Fixed points of twisted Coxeter element equals centre
%%%
%%%
\begin{corollary}
i) $\mbox{dim}\{x\in\csa , \mbox{cox}^{\tau }(x)=x\}=1$ and\\
ii) $\mbox{dim}\{x\in\csa' , \mbox{cox}^{\tau }(x)=x\}=1$.
\end{corollary}
Since we know that the line $\C c=\C\sum_{i=0}^r\check{a}_i\check{\alpha}_i$ is fixed by the $s_i$ as well as $\tau $ we must have:
\begin{equation}
\{x\in\csa , \mbox{cox}^{\tau }(x)=x\}=\{x\in\csa' , \mbox{cox}^{\tau }(x)=x\}=\C c.
\end{equation}
In other words, the twisted Coxeter element has eigenvalue one with multiplicity one on $\csa $
as well as on $\csa'$. 

Furthermore, since the element $\delta$ is fixed by the $s_i$ and by $\tau $ we have:
%%%
%%%
%%% Multiplicities of eigenvalue 1 of twisted Coxeter element
%%%
%%%
\begin{eqnarray}
\tau (d)&=&d+h, \mbox{with}\;h\in\csa'\\
s_i(d)&=&d+h_i,\mbox{with}\;h_i\in\csa'.  
\end{eqnarray}
Hence, $\mbox{cox}^{\tau }$ has upper triangular shape with respect to the direct sum decomposition
$\csa =\C d\bigoplus \csa'$.
Thus, $1$ is at least a zero of the characteristic polynomial of $\mbox{cox}^{\tau }$ on $\csa $ with
multiplicity $2$.
In fact, even more is true:
%%%
%%%
%%% Proposition 2.2: Multiplicity of principal spaces of twisted Coxeter element
%%%
%%%
\begin{proposition}
The multiplicity of $1$ as zero of the characteristic polynomial of $\mbox{cox}^{\tau }$ on $\csa' $
is two and hence, is three on $\csa $.
\end{proposition}
%%%
%%%
%%% Characteristic polynomials of twisted Coxeter elements
%%%
%%%
For $\tau=id$ this statement was proven by Br\"uchert, \cite{bruechert} Proposition 9, by direct computation
using the characteristic polynomials which were computed by Coleman \cite{coleman:coxetertrafos}, Table 3.
For non-trivial $\tau$  it also follows by direct computation of the characteristic polynomials of the twisted Coxeter
elements for all cases. These calculations have been carried out by the author and the results are listed in Table 2 in the appendix.
%%%
%%%
%%%  A^1_n case
%%%
%%%
In the $\mbox{A}^1_n,n\geq 2$ we set $\mbox{cox}^{\tau }=s_{i_1}...s_{i_p}\tau$, where the $i_j$ have smallest 
index among all indices in its orbits and they are in increasing order $i_1<i_2<....<i_p$. 
(Recall that the twisted Coxeter element is not unique up to conjugation in this case.)

%%%
%%%
%%%  Solution to coxtw-1(b)=kc
%%%
%%%
For the construction of the $\C^{\ast }$-action on the Steinberg cross section 
in section 4 we need a solution $b\in\csa '$ to the equation $(\mbox{cox}^{\tau }-1)(b)=c
=\sum_{i=0}^r\check{a}_i\check{\alpha}_i$. 
This equation has clearly a solution over the space $\Q <\check{\alpha}_0,...
,\check{\alpha}_n>\subset \csa'$. We can even achieve $b\in\Z\check{\widehat{\Pi}}$ by multiplying the above equation with 
an appropriate integer.
%%%
%%%
%%%  Conventions for twisted Coxeter element
%%%
%%%
However this solution clearly depends on the explicit choice of $\mbox{cox}^{\tau }$.
Let us make the following choices for the twisted Coxeter elements:
In the $\mbox{A}^1_{n}$ case we use $\mbox{cox}^{\tau } $ as above. 
For the cases $\mbox{B}^1_{n}$, $\mbox{A}^2_{2n-1}$ we take $\mbox{cox}^{\tau }=
s_1...s_n\tau $. In the case of $\mbox{E}^1_6$ and the automorphism $\tau $, as above,
we choose $\mbox{cox}^{\tau }=s_1s_2s_3\tau $.  In the $\mbox{D}^1_{n} $ we make the following
choices: if $0$ and $1$ lie in the same orbit of the exterior automorphism we take representatives of 
$\Sigma $-orbit of minimal possible index in $I\backslash\{0\}$ and multiply them in 
the order of increasing index. Finally, we multiply the result with the automorphism $\tau $ from the right.
In all the other cases, ie. $\mbox{C}^1_{n}$, $\mbox{D}^2_{n+1}$, $\mbox{E}^1_7$, $\mbox{E}^1_6$ 
and $\sigma $ and all the other automorphisms of $\mbox{D}^1_{n} $ we choose representatives
of the $\Sigma $ orbits of least possible index which multiply in increasing order and finally we 
multiply the result with the automorphism $\tau $ from the right.
Then by direct computation we obtain the following generalisation of Proposition 10 in \cite{bruechert} for $\tau \neq id$:
%%%
%%%
%%%  Proposition 2.3: Description of b
%%%
%%%
\begin{proposition}
\label{descriptionb}
Let $\widehat{\Pi}$ be an affine Dynkin diagram,  $\widehat{\W }$ its Weyl group, and $\tau\in\mbox{Aut}(\widehat{\Pi}) $ 
an exterior automorphism. Then, for the twisted Coxeter element $\mbox{cox}^{\tau }$ the following
holds:
There is an element $b\in\N<\check{\widehat{\Pi} }>$ such that $(\mbox{cox}^{\tau }-1)(b)=kc$ which is uniquely
determined by requiring $k\in\N $ to be minimal and $b-c\notin \N<\check{\widehat{\Pi }}>$.
The number $k$ as well as the element $b$ expressed in the basis $\check{\widehat{\Pi}} $ 
are listed in the Table 3 in the appendix for $\tau\neq id$.
\end{proposition}
%%%
%%%
%%%  The betas evaluated on b
%%%
%%%
Now, consider those positive roots $\beta_0,...,\beta_{s-1}$ mapped to negative ones by $\mbox{cox}^{\tau \,-1}=\tau^{-1}s_{s-1}...s_{0} $. 
It is easy to see that these roots can be written as: $\beta_i=s_0...s_{i-1}\alpha_{i}$, $i\in\{0,...,s-1\}$.
Then, the following result holds:
%%%
%%%
%%% Lemma 2.4: betas on b are proportional to folded dual Kac labels
%%%
%%%
\begin{lemma}
\label{twistedKaclabel}
The $s$-tuple $(\beta_0(b),...,\beta_{s-1}(b))$ coincides up to a factor with the $s$-tuple of the dual Coxeter numbers 
of the folded Dynkin diagram ${}^{\tau}\widehat{\Pi}$. (In the cases where there is no folded Dynkin diagram, ie. the cases 
$\mbox{A}_n^1,\tau$ with $\tau =\gamma$ generating $\widehat{P}/\widehat{Q}$ we set $a^{\tau}_0=1$.) 
\end{lemma}
%%%
%%%
%%%  Proof of Lemma 2.4
%%%
%%%
{\it Proof:} First, consider the following statement which follows by an easy induction:

{\sl Claim:} Let $w\in \widehat{\W}$ an element with reduced expression $w=s_{i_1}...s_{i_t}$. Then for any $\lambda\in\csa^{\ast}$
we have the equation:
\begin{equation}
\lambda-w(\lambda)=\sum_{j=1}^{t}\lambda(\check{\alpha}_{i_j})\,s_{i_1}...s_{i_{j-1}}(\alpha_{i_j}).
\end{equation}
Now, set $\tilde{\alpha}_i:=\frac{1}{|\Sigma_i|}\sum_{l=0}^{N-1}\tau^{l}(\alpha_i)$.
Applying the formula above to the case $w=s_{0}...s_{s-1}\in\widehat{\W}$ and 
$\lambda=\tau^{l}(\alpha_i)$ yields:
\begin{eqnarray}
\tilde{\alpha}_i-\mbox{cox}^{\tau}(\tilde{\alpha}_i)&=&\tilde{\alpha}_i-s_0...s_{s-1}(\tilde{\alpha}_i)=\\
&=&\frac{1}{|\Sigma_i|}\sum_{l=0}^{N-1}(\tau^{l}(\alpha_i)-s_0...s_{s-1}(\tau^{l}(\alpha_i))=\\
&=&\frac{1}{|\Sigma_i|}\sum_{l=0}^{N-1}\sum_{j=0}^{s-1}\widehat{C}_{j\tau^l(i)}\beta_j.
\end{eqnarray}
Evaluating the result on the element $b$ yields:
\begin{eqnarray}
\frac{1}{|\Sigma_i|}\sum_{l=0}^{N-1}\sum_{j=0}^{s-1}\widehat{C}_{j\tau^l(i)}\beta_j(b)=\tilde{\alpha}_i(b)-\mbox{cox}^{\tau}(\tilde{\alpha}_i)(b)=\\
\tilde{\alpha}_i(b-\mbox{cox}^{\tau\,-1}(b))=\tilde{\alpha}_i(kc)=0.
\end{eqnarray}
Here, the second last equation follows from the definition of $b$ using the fact that $\mbox{cox}^{\tau }$ fixes $c$.
Thus the $s$-tuple $(\beta_0(b),...,\beta_{s-1}(b))$ is mapped to zero by the matrix $\mbox{diag}(t_0^{-1},...,t_{s-1}^{-1})\:{}^{\tau}\widehat{C}$.
Recall that the $t_i$ are the numbers used in Equation \ref{eqnfoldeddiagram}. 
This clearly implies the result. \hfill $\qed $.

Using the list of the dual Kac labels we can calculate explicitly the factor $p$ such that 
$(\beta_0(b),...,\beta_{s-1}(b))=p(\check{a}_0^{\tau },...,\check{a}_{s-1}^{\tau })$. 
These values are compiled in Table 4 of the appendix.
\section{Affine Kac-Moody Lie Algebras and Groups}
%%%
%%%
%%%  Introductory part  
%%%
%%%
In this section we recall the basic notions of the theory of affine Kac-Moody Lie algebras and Groups and their representation theory.
However, this section contains new results on central extensions of non-connected loop groups.
%%%
%%%
%%%  Subsection 3.1: Affine lie algebras.
%%%
%%%
\subsection{Affine Lie algebras}
Let $\liea $ be a complex finite dimensional simple Lie algebra  of rank $r$ and $\csa\subset\liea $ its Cartan 
subalgebra. We denote by $\Delta $ its root system and by $\Pi:=\{\alpha_1,...,\alpha_r\} $ the set of simple 
roots.
 Then, ${\cal L}\liea :=\{x:\C^*\to \liea, x\; holomorphic\}$ is the (holomorphic) loop algebra. 
Now, it is a well known fact that all central extensions $\widetilde{{\cal L}\liea}^k $, $k\in \C $ of ${\cal L}\liea $ by 
$\C$ can be described by the following commutation relations: 
\begin{eqnarray}
 \left[c,x\right]     &=&0 \\             
 \left[x,y\right](z) &=&\left[x(z),y(z)\right]+\frac{k}{2\pi i} \int _{|z|=1}<\frac{d}{dz}x(z),y(z)>dz \;c
\end{eqnarray} 
Here, $c$ denotes the central variable of $\widetilde{{\cal L}\liea}^k $ and $<,>$ is the Killing form normalised in 
such a way that $<\alpha,\alpha>=2$ for all long roots $\alpha\in\Delta$.
Let us mention as an aside that for $k\neq 0$ all the $ \widetilde{{\cal L}\liea}^k $ are isomorphic as Lie 
algebras while they are not as central extensions.
Observe that the derivation $d=z\frac{d}{dz}$ on $\looplie $ lifts to a derivation on $ \widetilde{{\cal L}\liea}^k $
by acting trivially on the centre $\C c $. The corresponding semidirect product $\widetilde{{\cal L}\liea}^k\oplus \C d $
is called {\it non-twisted affine Kac-Moody Lie algebra } and denoted by $\widehat{{\cal L}\liea}^k$.
If we start with polynomial loops instead of holomorphic ones a similar construction yields a polynomial version
$\kmliea_{\tiny\mbox{pol}}=\liea\otimes\C[z,z^{-1}]\oplus \C c\oplus \C d$ of the algebra above.
This is the non-twisted affine Lie algebra  in the sense of Kac \cite{kac:infdimliealgebras}. It was pointed out in 
\cite{goodmanwallach:unitarystructure} that 
$ \kmliea $ is a completion of $\kmliea_{\tiny\mbox{pol}} $ in a certain topology. 
Similar to the finite dimensional case the polynomial version has a nice root space decomposition: 
Note that the subalgebra $\widehat{\csa }:=\csa \oplus \C c\oplus \C d $ is a Cartan subalgebra of $\kmliea_{\tiny\mbox{pol}} $. 
This allows us to define a root system $\widehat{\Delta}$ which is an affine root system. 

The corresponding root space decomposition looks as follows:
\begin{equation}
\kmliea_{\tiny\mbox{pol}} = \widehat{\csa }\oplus\bigoplus_{\tilde{\alpha }\in \widehat{\Delta}}\kmliea_{\tilde{\alpha }}.
\end{equation}
Here, the root spaces are: 
$\kmliea_{\tilde{\alpha }}=\liea_{\alpha }\otimes z^n $, for $\tilde{\alpha }=\alpha + n\delta $ and $\kmliea_{\tilde{\alpha }}=
\csa\otimes z^n$, for $\tilde{\alpha }=n\delta $ with $\alpha\in\Delta$.

Let us give a few remarks about twisted affine Lie algebras. Let $\sigma $ be an automorphism of
$\liea $ of order $s$. Then consider the fix point subalgebra 
\begin{equation}
 \looplie(\sigma):=\{x\in\looplie, \sigma(x(e^{-\frac{2\pi i}{s}}z)=x(z)\}.
\end{equation}
Using a similar construction as above yields the {\it twisted affine Kac-Moody Lie  algebra } $\kmliea\!\! (\sigma )$.
As in the non-twisted case we have the notion of root systems, Kac labels, etc. For the details we refer to 
\cite{kac:infdimliealgebras}.

We conclude this section by the following definition, already introduced in \cite{fuchs}:
\begin{definition}
Let $\kmliea $ be any twisted or untwisted affine Kac-Moody Lie algebra with Dynkin diagram $\widehat{\Pi}$ and let $\tau\in\mbox{Aut}(\Pi) $ 
be a automorphism thereof.
Then the Lie algebra $\widehat{{\cal L}\liea }^1({}^{\tau}\widehat{\Pi})$ corresponding to the folded Dynkin diagram is called the orbit 
Lie algebra of the pair $(\widehat{\Pi},\,\tau)$.
\end{definition}
{\it Remark:} (i) Note, that the orbit Lie algebra can in general not be realised as a subalgebra of $\kmliea $, not even for $k=1$.

(ii) The orbit Lie algebra will play an important role in the representation theory, see Subsection 3.3.
%%%
%%%
%%%  Subsection 3.2: Affine Kac Moody Groups
%%%
%%%
\subsection{Affine Kac Moody Groups}
In this subsection we will introduce affine Kac Moody groups along the long line as in the Lie algebra case. 
Let $G$ be an algebraic group over $\C$ with simple identity component $G_0$ and
Lie algebra $\liea $. If $G$ is not connected itself we assume that $G$ is a semidirect product $G_0\rtimes\Gamma $
with $\Gamma<\mbox{Aut}(\Pi)$ being a subgroup of the automorphism group of the Dynkin diagram $\Pi $ of $G_0$.
Denote by $\tilde{G}$ the universal cover of $G_0$. 
Then, $\pi_1(G_0)$ is the fundamental group of $G_0$.
There is a central subgroup $ Z\subset C(G_0)\subset \tilde{G}$ such that $G_0=\tilde{G}/Z$ and $G=\tilde{G}\rtimes\Gamma/Z$ .
We define the {\it loop group } $\loopgr $ by $\loopgr :=\{X:\C^*\to G, X holomorphic \}$.
Endowed with the compact-open topology it becomes a topological group. 
On $\loopgr $ we have an action of $\C^*$ by rotations:
\begin{eqnarray}
\C^*\times \loopgr &\to & \loopgr\\
(s,X)&\mapsto & s.X,\mbox{ with:}\, (s.X)(z)=X(s^{-1}z).
\end{eqnarray}
It is well known fact that the component group $\pi_0(\loopgr ) $ of $\loopgr $ can be described by the equation:
\begin{equation}
\pi_0(\loopgr )=\pi_1(G_0)\rtimes\Gamma.
\end{equation}
Embedding $G\to\loopgr $
as the constant loops provides us with an isomorphism ${\cal L}\tilde{G}/Z\cong \loopgr_0$ where the right hand side 
is the identity component of $\loopgr $. As will be shown below the loop group itself is a semidirect product:
\begin{equation}
\loopgr =\loopgr_0\rtimes \pi_0(\loopgr )= \loopgr_0\rtimes (\pi_1(G_0)\rtimes\Gamma ).
\end{equation}
Since the central extensions of non-connected loop groups are rather complicated to handle we will postpone this and describe briefly the 
connected case first. For details we refer to \cite{presslaysegal:loopgroups}. 
Let $G$ be connected and simply connected, ie. $G=\tilde{G}$. Every central extensions of $\loopgr $ by $\C^*$ gives rise to a central 
extensions of its Lie algebra $\looplie $ by $\C$. But, a central extension of $\looplie $ only lifts to the group if  
the integrality condition $k\in \Z$ for our parameter $k$ as above, is satisfied.
The corresponding central extension of $\loopgr $ is denoted by $\widetilde{{\cal L}G}^k$.
One observes that the rotation action of $\C^*$ on $\loopgr $ lifts to $\widetilde{{\cal L}G}^k$. 
The semidirect product $\widetilde{{\cal L}G}^k\rtimes\C^* $ is called 
{\it affine Kac-Moody group of level $k$} and denoted by $\kmG $.

%%%
%%%
%%%  The non simply connected case
%%%
%%% 
Now, consider the case of not necessarily simply connected and connected $G$ with component group $\Gamma $ as above. 
Our discussion will follow the line of reasoning suggested by Toledano Laredo \cite{toledanolaredo:representations}
where the connected but non simply connected case is treated. However, for the construction of central extensions we 
will restrict ourselves to a cyclic subgroup of the loop group's component group.

Instead of working with the loop group itself we introduce the group of `open loops' with `boundary values' in $Z$ :
\begin{equation}
\goloopgr :\{X:\C\to \tilde{G}\rtimes\Gamma, X\mbox{holomorphic}, X(t)X(t+1)^{-1}\in Z\}.
\end{equation}
Let us start by discussing its structure:
This group is clearly a semidirect product $\oloopgr\rtimes\Gamma $ with 
$\oloopgr :=\{X:\C\to \tilde{G}, X\mbox{holomorphic}, X(t)X(t+1)^{-1}\in Z\}$ where the $\Gamma $-action is defined by 
$\sigma(X)(t)=\sigma(X(t))$ for $\sigma\in\Gamma $.
Furthermore, we easily see that $\pi_0(\oloopgr )=Z$ and that $\goloopgr_0={\cal L}\tilde{G}$. 
Identifying $z$ with $e^{2\pi i t}$ yields an isomorphism $\loopgr\cong\oloopgr /Z$.
But, the group $\oloopgr $ itself is a semidirect product as well.
Indeed, fix a maximal torus $\tilde{T}$ of $\tilde{G}$ which is mapped to a maximal torus $T:=\tilde{T}/Z$ of $G$
and denote the corresponding co-character lattices by $ \check{\chi}(\tilde{T})$ respectively $\check{\chi}(T)$, 
ie. $\check{\chi}(\tilde{T}):=Hom_{\mbox{\tiny alg gr.}}(\C^*, \tilde{T})$.
We can identify every element $\check{\beta}\in\check{\chi}(T) $ with the open loop,
 also denoted by $\check{\beta}$: 
\begin{equation}
\check{\beta}(t):=exp (2\pi i t\check{\beta }).
\end{equation} 
Here we regard $\check{\beta}$ as an element in the real part $\csa_{\R}=\R \check{\Pi}$ of the Cartan subalgebra, and  
$exp$ is the exponential map of $\tilde{G}$.  
By this identification, we obtain $\check{Q}=\check{\chi}(\tilde{T})\subset {\cal L}\tilde{G}$. 
Thus, the following lemma holds:
\begin{lemma}
The group of open loops has the following semidirect product structure:
\begin{eqnarray}
\oloopgr\cong({\cal L}\tilde{G}\rtimes \check{\chi}(T))/\check{\chi}(\tilde{T})\cong{\cal L}\tilde{G}\rtimes Z.\\
\goloopgr \cong ({\cal L}\tilde{G}\rtimes Z)\rtimes\Gamma \cong {\cal L}\tilde{G}\rtimes (Z\rtimes\Gamma ).
\end{eqnarray}
\end{lemma}
Since we are only interested in a single exterior component of $\goloopgr $
we restrict ourselves to a cyclic subgroup $\Sigma \subset\pi_0(\goloopgr )=Z\rtimes\Gamma $.
Assume that $\Sigma $ is generated by $\tau =\rho\sigma $ with $\rho\in Z $ and $\sigma\in\Gamma $.
 
Our next goal is to construct the central extensions of ${\cal L}\tilde{G}\rtimes\Sigma\subset
\goloopgr $: 
Using the identification $Z=\check{\chi}(T)/\check{\chi}(\tilde{T})$, there is a fundamental co-root
$\check{\lambda}_{\rho}$ (respectively $\check{\lambda}_{\rho}=0$, for $\rho=id $) 
such that $\check{\lambda}_{\rho}+\check{\chi}(\tilde{T})=\rho $.
Considering $\check{\lambda}_{\rho}$ as an open loop as introduced above the element $\check{\lambda}_{\rho}\sigma {\cal L}\tilde{G}$
generates the component group of ${\cal L}\tilde{G}\rtimes\Sigma $ which is $\Sigma $.
Indeed, if $p$ is the order of the diagram automorphism $\rho\sigma $ 
we see that $(\check{\lambda}_{\rho}\sigma )^p\in\chi(\tilde{T})$.
Furthermore, we have $1=(\rho\sigma )^p=\rho\sigma(\rho)...\sigma^{p-1}(\rho)\sigma^p$
which implies $\sigma^p=1$. Thus, $(\check{\lambda}_{\rho}\sigma )^p=
\check{\lambda}_{\rho}+\sigma(\check{\lambda}_{\rho})+...+\sigma^{p-1}(\check{\lambda}_{\rho})$
has even to be $\sigma $-invariant.

The construction of the central extension proceeds in several steps:
First, recall that all the isomorphism classes of central extensions of any group $H$ by an abelian group $A$ are given by the 
second cohomology group $H^2(H,A)$. 
Now, observe that the exact sequence
\begin{equation}
1\longrightarrow {\cal L}\tilde{G} \longrightarrow  {\cal L}\tilde{G}\rtimes \Sigma\longrightarrow\Sigma \longrightarrow 1.
\end{equation}
gives rise to the following long exact cohomology sequence (here, we use the fact 
$\mbox{Hom}({\cal L}\tilde{G_0},\Cst)=\{id\}$):
\begin{equation}
1\longrightarrow H^2(\Sigma,\Cst ) \longrightarrow H^2({\cal L}\tilde{G}\rtimes \Sigma,\Cst )\longrightarrow 
H^2({\cal L}\tilde{G_0},\Cst )\longrightarrow ...
\end{equation} 
Since $\Sigma $ is cyclic we have $H^2(\Sigma,\Cst)=0$.  By the classification of its central extensions,
as mentioned in the beginning of this section, $H^2({\cal L}\tilde{G},\Cst)=\Z$. 
Thus, $H^2({\cal L}\tilde{G}\rtimes \Sigma,\Cst )$ is either isomorphic to $\Z$ or trivial.
Our next goal is to show that it is not trivial. 

This follows from the explicit construction of its central extensions: We begin by  
lifting the action of the subgroup $\tilde{\Sigma }:=<\check{\lambda}_{\rho}\sigma >$ to the centrally extended 
loop group $\widetilde{ {\cal L}\tilde{G}}^k$. Note, that the group $\tilde{\Sigma }$ is $\Z$ unless $\rho=id $ in which case it is 
finite.
By direct computation and Lemma 3.1.1 in \cite{toledanolaredo:representations} we obtain:
\begin{lemma}
The conjugacy action of $\check{\lambda}_{\rho}\sigma $ on $\looplie $ resp. ${\cal L}\tilde{G}$ lifts uniquely to 
any of their central extension. For any $x+ac\in\widetilde{{\cal L}\liea}^k $ it is given by:
\begin{equation}
\label{actiontau}
\widetilde{\mbox{Ad}\check{\lambda}_{\rho}\sigma }(x+ac)=\mbox{Ad}(\check{\lambda}_{\rho}\sigma )(x)+(a+\frac{k}{2\pi i}
\int_{t=0}^1<\check{\lambda}_{\rho}^{-1}\frac{d}{dt}\check{\lambda}_{\rho},\sigma(x)>dt )c.
\end{equation}   
In particular, for $\rho=id $ this means that the integral vanishes. 
\end{lemma}
This Lemma allows us to form the semidirect product 
$\widetilde{{\cal L}\tilde{G}}^k\rtimes \tilde{\Sigma }$ for the central extension of 
${\cal L}\tilde{G}$ of level $k$.
The $\widetilde{\Sigma}$-action is the lift of its action on the Lie algebra, as given in Equation \ref{actiontau}.

Now, consider the $ \sigma$-invariant co-root $\check{\beta}:=(\check{\lambda}_{\rho}\sigma)^p$
and denote by $\tilde{ \check{\beta}}\in \widetilde{{\cal L}\tilde{G_0}}$ an arbitrary lift of it. 
A calculation as in the proof of \cite{toledanolaredo:representations}, Proposition 3.2.2 yields:
\begin{equation}
\widetilde{\mbox{Ad}(\check{\lambda}_{\rho}\sigma)}(\tilde{ \check{\beta}})=(-1)^{k<\check{\lambda}_{\rho},\check{\beta}>}.
\label{centralextension}
\end{equation}

Observe that the set $\widetilde{\Sigma}^p=\{\check{\beta}^n,\,n\in\Z\}\subset {\cal L}\tilde{G}$ is isomorphic to $\Z$. 
Since $H^2(\Z,\Cst)$ is trivial $\widetilde{{\cal L}\tilde{G}}^k|_{\widetilde{\Sigma}^p}$ has to be trivial as well. Thus $\widetilde{\Sigma}^p$
can be embedded into $\widetilde{{\cal L}\tilde{G}}^k$. 

For the existence of a central extension of ${\cal L}\tilde{G}\rtimes\Sigma $ it is necessary that the two elements $\check{\lambda}_{\rho}\sigma $
and $\check{\beta}$ should commute, ie. the following condition should hold:

\begin{equation}
(-1)^{k<\check{\lambda}_{\rho},\check{\beta}>}=1.
\label{compitibilycond}
\end{equation}

Following the reasoning of \cite{toledanolaredo:representations}, proof of Proposition 3.3.1, Equation \ref{compitibilycond} implies
that the subgroup $N:=\{(\check{\beta}^n,\,(\check{\lambda}_{\rho}\sigma )^{-pn},\,n\in\Z\}$ is a normal subgroup of 
$\widetilde{{\cal L}\tilde{G}}\rtimes \widetilde{\Sigma} $. Therefore, $N$ is normal
if and only if $k<\check{\lambda}_{\rho},\check{\beta}>\in 2\Z $.
In this situation we define the central extension of level $k$ of ${\cal L}\tilde{G}\rtimes \Sigma $ by:
\begin{equation}
\label{centralext}
\widetilde{{\cal L}\tilde{G}\rtimes \Sigma}^k:=(\widetilde{{\cal L}\tilde{G}}^k\rtimes \widetilde{\Sigma})/N
\end{equation}

\begin{definition} 
The smallest positive integer $k$ such that 
$k<\check{\lambda}_{\rho},\check{\beta}>\in 2\Z $ is called the
fundamental level of ${\cal L}\tilde{G_0}\rtimes \Sigma $ and will be denoted by $k_f$.
\end{definition}

Cohomologically expressed this translates to the following exact sequence:
\begin{equation}
1 \longrightarrow H^2({\cal L}\tilde{G}\rtimes \Sigma,\Cst )\longrightarrow 
H^2({\cal L}\tilde{G},\Cst )=\Z\longrightarrow \Z/k_f\Z\longrightarrow 1.
\end{equation}

We summarise the classification of central extensions:

\begin{theorem}
\label{classoscentralexts}
Every central extension of  ${\cal L}\tilde{G}\rtimes \Sigma $ is uniquely determined by
the level $k$ of its restriction to ${\cal L}\tilde{G}$ which has to be 
a multiple of $k_f$.\\
The corresponding central extension is described by Equation \ref{centralext}.
\end{theorem}

{\it Remarks:} i) Let us illustrate the condition $k_f|k$ of allowed level by the following example: 
Consider $G=PSl_2$, ie. $\tilde{G}=Sl_2$ and 
$\Sigma=Z=\{\pm id\}$. Then, we have $\widetilde{\Sigma}=\Z\frac{\alpha}{2}$, $p=2$, 
$\check{\lambda}_{\rho}=\frac{\check{\alpha}}{2}$ and thus $\widetilde{\Sigma}^2=\Z\alpha $.
Now Equation \ref{centralextension} yields $(-1)^{k<\check{\alpha},\frac{\check{\alpha}}{2}>}=-1$
for $k$ odd violating condition \ref{compitibilycond}.
Thus, there are no central extensions of ${\cal L}_{\Z/2\Z}Sl_2$ of odd level.

ii) It is easy to see that $k_f$ is either $1$ or $2$.

iii) If $\rho= id $ then $k_f=1$.

iv) In \cite{toledanolaredo:representations} the classification of central extensions 
of $\oloopgr $ is given by Proposition 3.3.1. In his result there appears an additional skew symmetric
form $\omega $ on $\check{\chi}(T)$ with values in $\Cst $. This form classifies 
central extensions of $\check{\chi}(T)$. Furthermore, a compatibility requirement between 
$\omega $ and the level $k$ appears, looking similar to Equation \ref{compitibilycond}.
However, it turns out that for cyclic $Z$ this form $\omega $ is uniquely determined by the level $k$. 
Therefore, we do not need $\omega $ in our approach.

v) the fundamental levels were calculated in \cite{toledanolaredo:representations}. For the reader's 
convenience we will summarise them in Table 5 in the Appendix.

The final step in the construction of the Kac-Moody group is to lift the $\Cst$-action to the central extension.
This will be carried out along the lines of Section 3.4 in \cite{toledanolaredo:representations}.
First, recall that $w\in\Cst $ acts on $X\in\loopgr $ by $(w.X)(z)=X(w^{-1}z)$. Now, observe that $\Cst $  does no longer act on 
$\goloopgr $. However, its universal cover $\C $ does:  
\begin{equation}
(s.X)(t):=X(t-s),\quad\forall s,t\in\C\,\mbox{and}\,X\in\goloopgr.
\end{equation}
Observe that this action restricted to ${\cal L}\tilde{G}$ factors through the usual
$\Cst=\C/\Z$-action. Using the fact that $\sigma $ commutes with the $\C$-action
we obtain for the action of $s\in\C $ on $\check{\lambda}_{\rho}\sigma$ the formula:
\begin{equation}\label{translationeqn}
s.\check{\lambda}_{\rho}\sigma=exp(-2\pi i s\check{\lambda}_{\rho})\, \check{\lambda}_{\rho}\sigma.
\end{equation}
This allows us to rewrite the semidirect product decomposition $({\cal L}\tilde{G_0}\rtimes \Sigma)\rtimes\C$ in the following way:
\begin{equation}
({\cal L}\tilde{G}\rtimes \Sigma)\rtimes\C= ({\cal L}\tilde{G}\rtimes \C)\rtimes\Sigma .
\end{equation}

Using \cite{toledanolaredo:representations} Lemma 3.1.1. we see that any central extension of ${\cal L}\tilde{G}\rtimes \C$ 
is uniquely determined by its restriction $\widetilde{{\cal L}\tilde{G_0}}$ to ${\cal L}\tilde{G}$, 
(up to a character of $\C$ which we fix to be the identity). Furthermore, the central extension is given by a semidirect product
$\widetilde{{\cal L}\tilde{G}\rtimes \C}=\widetilde{{\cal L}\tilde{G}}\rtimes \C$. Now, let $k$ with $k_f|k$ be the level of this 
central extension.

For lifting the action of $\check{\lambda}_{\rho}\sigma$ to this centrally extended group
we consider the action on its Lie algebra
$\mbox{Lie}\widetilde{{\cal L}\tilde{G}\rtimes \C}^k=\widetilde{\looplie}^k\oplus \C d$. 
Here $d$ is given by $\frac{1}{2\pi i}\frac{d}{dt}=z\frac{d}{dz}$.

A simple calculation using the fact that $\sigma$ fixes $d$ yields:
\begin{eqnarray}
\widetilde{\mbox{Ad}(\check{\lambda}_{\rho}\sigma )}(x+bd+ac)&=& \mbox{Ad}\check{\lambda}_{\rho}(x)+\frac{b}{2\pi i}\check{\lambda}_{\rho}
\frac{d}{dt}\check{\lambda}_{\rho}^{-1}+bd+\\
&+&(a+\frac{k}{2\pi i}\int_{t=0}^1<\check{\lambda}_{\rho}^{-1}\frac{d}{dt}\check{\lambda}_{\rho},\,\sigma(x)>dt\notag\\ 
&+&b\frac{k}{8\pi^{2} }\int_{t=0}^1<\check{\lambda}_{\rho}^{-1}\frac{d}{dt}\check{\lambda}_{\rho}\,
\check{\lambda}_{\rho}^{-1}\frac{d}{dt}\check{\lambda}_{\rho}>dt)c.\notag
\end{eqnarray}

This allows us to construct the semidirect product $\widetilde{{\cal L}\tilde{G}\rtimes \C}^k\rtimes\widetilde{\Sigma}$.

Similar to the situation discussed above we deduce that the group 
$N:=\{(\check{\beta}^n,1,\check{\beta}^{-n}),\,n\in\Z\}\subset(\widetilde{{\cal L}\tilde{G}}^k\rtimes \C)\rtimes \widetilde{\Sigma}$ 
is a normal subgroup.

\begin{definition}
The quotient group $\widehat{{\cal L}_{\Sigma}\tilde{G}}^k:=((\widetilde{{\cal L}\tilde{G}}^k\rtimes \C)
\rtimes \widetilde{\Sigma})/N$ is called the affine Kac-Moody group corresponding to $G$ and $\Sigma$.
\end{definition}
{\it Remark:}
Let $q$ be the smallest positive integer such that $q\check{\lambda}_{\rho}\in\check{Q}$.
It follows from Equation \ref{translationeqn}
that $\check{\lambda}_{\rho}\sigma $
fixes the translation by $q$. Thus we, effectively,  obtain an action of $\widetilde{\Cst }:=\C/q\Z $.

Let us briefly discuss the case of twisted Kac-Moody groups. Consider a simple (and thus connected)
algebraic group $G$ and $\sigma $ one of its exterior
automorphisms of order $s$. The twisted loop group $\loopgr (\sigma)$ is given as fixed point group in the following way:
\begin{equation}
\loopgr (\sigma):=\{X\in\loopgr | \sigma(X(e^{\frac{2\pi i}{s}}z))=X(z)\}.
\end{equation} 
This group is known to have $|\pi_1(G)^{\sigma}|$ many connected components.
Thus, the only non connected twisted loop groups may appear for $G$ of case $\mbox{A}_{2n+1}$ and $\mbox{D}_{n+1}$ yielding that 
the number of connected components is at most $2$. This reflects the fact that the only twisted affine Dynkin diagrams
with nontrivial symmetries are of this type. Denote the nontrivial symmetry by $\gamma $.

As in the non-twisted case let us first review the case of simply connected $G$: Also in this situation the central
extensions of $\loopgr (\sigma)$ are classified by their level $k\in\Z$. We denote them by $\widetilde{\loopgr (\sigma)}^k$.
Furthermore, the restriction of the $\Cst $ action to $\loopgr (\sigma)$ can be lifted to these central 
extension. Then the Kac-Moody group $\widehat{\loopgr (\sigma)}^k$ is given by:
\begin{equation}
\widehat{\loopgr (\sigma)}^k:=\widetilde{\loopgr (\sigma)}^k\rtimes\Cst.
\end{equation} 

If $\loopgr (\sigma)$ has two connected components our discussion will follow the line of reasoning of the
non-twisted case. Therefore, we shall give only a brief account:
Let $\tilde{G}$ be universal cover of $G$ and let $Z\subset\tilde{G}$ its central subgroup such that $G=\tilde{G}/Z$ . 
As above choose maximal tori $\tilde{T}$ of $\tilde{G}$ and $T$ of $G$.
Denote by $\check{\chi}(T)^{\sigma}$ respectively $\check{\chi}(\tilde{T})^{\sigma}=\check{Q}^{\sigma} $
the corresponding fixed point character lattices. 
Introduce the twisted `group of open loops' with `boundary values' in $Z^{\sigma }$:

\begin{equation}
\oloopgr (\sigma ):=\{X\in {\cal L}\tilde{G}, \sigma (X(e^{\frac{2\pi i}{s}}z))X(z)^{-1}\in Z^{\sigma}\}.
\end{equation} 

Let $\check{\lambda}_{\gamma}\in\check{\chi}(T)^{\sigma}$ be the fundamental co-weight corresponding to the non-trivial element of
$\gamma\in Z^{\sigma}$. Then, we have $\check{\alpha}_{\gamma}:=2\check{\lambda}_{\gamma}\in  \check{Q}^{\sigma}$. The loop 
$\check{\alpha}_{\gamma}:z\mapsto\check{\alpha}_{\gamma}(z)$ is contained in the exterior component of $\oloopgr (\sigma )$.
Indeed, consider $\sigma(\check{\alpha}_{\gamma}(e^{\frac{2\pi i}{2}}z))\check{\alpha}_{\gamma}(z)^{-1}=\check{\alpha}_{\gamma}(-1)$,
which is the non-trivial element in $Z^{\sigma}$.

There is again a semidirect product decomposition with $\Sigma:=\{id,\,\gamma\}$:
\begin{equation}
\oloopgr (\sigma ):=\loopgr(\sigma )\rtimes\Sigma.
\end{equation} 

Now, the non-connected twisted Kac-Moody group $\widehat{{\cal L}_{Z^{\sigma}}\tilde{G}}^k\!\!(\sigma)$
is constructed in the same way as in the non-twisted case with
$\widetilde{\Sigma }:=<\check{\alpha}_{\gamma}>$ and $p=2$.
Observe, that in this case the fundamental level is always $1$, since 
$<\check{\alpha}_{\gamma},2\check{\alpha}_{\gamma}>\in2\Z$.
%
%Our first goal is to construct all central extensions of $\oloopgr $ and then to compare these with the 
%pull backs of central extensions of $\loopgr $.
%
%It will turn out that central extensions  on $\oloopgr $ are given by central extensions of
%$\loopgr $ and $\chi(T) $ which coincide on $\chi(\tilde{T})$: 
%Any central extension $\widetilde{\chi(T)}$ of $\chi(T)$ by $\C^*$ is uniquely determined by 
%a skew symmetric $\Z$-bilinear form $\omega $, the so called {\it commutator map},
% on $ \chi(T)$ which is defined by:
%\begin{equation}
%\omega (\lambda, \, \mu )=\widehat{\lambda}\widehat{\mu}\widehat{\lambda}^{-1}\widehat{\mu}^{-1}.
%\end{equation}
%The elements $\widehat{\lambda}$ and $\widehat{\mu}$ are lifts of $\lambda, \,\mu\in\chi(T) $
%to $\widetilde{\chi(T)}$.
%Choose a central extension ${\widetilde{\cal L}^k\tilde{G}}$ of ${\cal L}\tilde{G}$ of level $k$.
%Assume that the commutator map of a chosen central extension $\widetilde{\chi(T)}$ of $\chi(T)$
%fulfils the condition:
%\begin{equation}
%\label{compatiblitycon}
%\omega (\lambda, \, \mu )=(-1)^{k<\lambda, \, \mu>}, \forall\lambda\in \chi(\tilde{T}), \mu\in\chi(T).
%\end{equation}

{\it Remark:} Before the conclusion of this subsection let us justify the use of the group of `open loops'.
We use the group of `open loops' and its corresponding Kac-Moody groups solely because they have the the right representation
theory. Toledano Laredo \cite{toledanolaredo:representations} indicates also a classification of 
all central extensions of the loop group ${\cal L}G$ with $G$ connected but not necessarily simply connected.
He showed that the pull back of any central extension of ${\cal L}G$ to ${\cal L}_ZG$ has to be of a level which is a multiple
of the basic level $k_b$. This basic level $k_b$ is a multiple of the fundamental level $k_f$. In our setting, ie. $Z=<\rho>$ cyclic, the basic level 
$k_b$ can be defined to be the smallest integer, such that $k_b<\check{\lambda}_{\rho},\check{\lambda}_{\rho}>$ is an integer.
If and only if this is the case the group $Z$ is contained in the centre of the centrally extended group of `open loops'.
The only further ambiguity entering in the classification of central extensions of ${\cal L}G$ is a character of $Z$.
Toledano Laredo showed that every projective highest weight representation of $\widetilde{{\cal L}_ZG}^{k_b}$
a level divisible by $k_b$ is already defined for $\widetilde{{\cal L}G}^{k_b}$. However, this is not true for 
the linear representation itself. Since we a working with the linear representations we have to use the group of
`open loops' as a kind of the ~universal cover'.
\subsection{Representation theory}
Here, we shall give a brief account on the representation theory of the Kac-Moody groups $\widehat{{\cal L}_{\Sigma}\tilde{G}}^k$
respectively $\widehat{{\cal L}_{Z^{\sigma}}\tilde{G}}^k\!\!(\sigma)$.
For the simplicity of the exposition we will formulate the results of this subsection only for the non-twisted case $\widehat{{\cal L}_{\Sigma}\tilde{G}}^k$.
However, analogues of all the stated results hold for $\widehat{{\cal L}_{Z^{\sigma}}\tilde{G}}^k\!\!(\sigma)$ as well.

The representation theory was already investigated in a couple of papers, see eg. \cite{fuchs,wendt:twistedcharacters,toledanolaredo:representations},
respectively \cite{jantzen:darstellungen} for the finite dimensional situation.

Irreducible highest weight representations of $\kmliea_{\scriptsize\mbox{pol}}$ are known to be labelled by the elements in $\csa^*$, their highest weight.
We denote such a representation by $V(\lambda)$ and call the number $\lambda(c)\in\C$ its level, see Kac's book \cite{kac:infdimliealgebras}.
Furthermore, the space $V(\lambda)_{\mu }$ will be the weight space of $V(\lambda)$ with respect to the weight $\mu\in\widehat{P}$.

For $\lambda\in\widehat{P}^+$ it was shown by Garland \cite{garland:loopalgebras} that $V(\lambda)$ carries a Hermitian from $(.,.)$ which is 
contravariant with respect to the Cartan involution on $\kmliea_{\scriptsize\mbox{pol}}$.
Note that the level is a non-negative integer in this case. From now on we assume that the level $k$ of the 
central extension divides the level $\lambda(c)$. 
Let us denote the Hilbert space completion of $V(\lambda)$ with respect to $(.,.)$ by $V(\lambda)^{ss}$.
Since the weight spaces $V(\lambda)_{\mu }$ are all finite dimensional $V(\lambda)^{ss}$ is the Hilbert space direct sum of them.
One observes that the action of $\kmliea_{\scriptsize\mbox{pol}}$ on $V(\lambda)$ does not extend to an action of
$\kmliea $ on
$ V(\lambda)^{ss}$. However, it was shown in \cite{goodmanwallach:unitarystructure,etingoffrenkelkirilov} that there is a 
dense subspace $V(\lambda)^{an}\subset V(\lambda)^{ss}$ containing $V(\lambda)$ carrying a $\kmliea $-action which extends the 
$\kmliea_{\scriptsize\mbox{pol}}$ on $ V(\lambda)$ by continuity.

Now, we turn to the representation theory of the affine Kac-Moody groups.
Let us first consider the simply connected case, ie. $G=\tilde{G}$. Then the following result is due to Etingof, Frenkel and Kirilov,
see \cite{etingoffrenkelkirilov}, Theorem 2.2 and Lemma 2.3:
\begin{theorem} 
(i) The action of the Lie algebra $\kmliea $ on $V(\lambda)^{an}$ integrates uniquely to an action of $\widehat{{\cal L}G}^k$.\\
(ii) For any $q\in\Cst$ with $|q|<1$ and any $g\in\widetilde{{\cal L}G}^k$ the operator $gq^{-d}:V(\lambda)^{an}\to V(\lambda)^{an}$
extends uniquely to a trace class operator on $V(\lambda)^{ss}$.
\end{theorem} 
We introduce the subset  $\widehat{{\cal L}G}^k_q:=\widetilde{{\cal L}G}^k\times\{q\}\subset \kmG $ and the submonoid
$\kmG_{<1}:=\cup_{q,\,|q|<1}\widetilde{{\cal L}G}^k_q\subset \kmG $. This subset is clearly invariant under conjugation. 
Furthermore, as a manifold it is isomorphic to $\widetilde{{\cal L}G}^k\times D^*$, where $D^*$ is the punctured unit disk in $\C$.
We will denote its elements by $(g,q)$ with $g\in\widetilde{{\cal L}G}^k $ and $q\in D^*$. By the theorem above every 
element $(g,q)\in\kmG_{<1}$ extends to a trace class operator on $V(\lambda)^{ss}$. This allows us to define for every 
$\lambda\in\widehat{P}^+$ its character $\chi_{\lambda}:\kmG_{<1}\to \C $ by the formula:
\begin{equation}
 \chi_{\lambda}(g,q):= Tr_{V(\lambda)^{ss}}(gq^{-d}).  
\end{equation}
These characters have the properties expected for characters, see \cite{etingoffrenkelkirilov}, Lemma 2.4 and Proposition 2.5:
\begin{proposition}
The functions $\chi_{\lambda}$ are holomorphic and invariant under conjugation.
\end{proposition}
Note that for the identification of the centre with $\Cst$, $\iota:\Cst\to \widetilde{{\cal L}G}^k$ we have 
$\chi_{\lambda}(\iota (u)g,q)=u^{\lambda(c)}\chi_{\lambda}(g,q)$.

Our next goal is to derive the representation theory for the non-connected Kac Moody groups. 
We start with a non-simply connected $G$ and consider the group $\widehat{{\cal L}_{\Sigma}\tilde{G}}^{k_f}$
as constructed in the preceding subsection for the fundamental level $k_f$. Let $V$ an irreducible 
module of $\widehat{{\cal L}_{\Sigma}\tilde{G}}^{k_f}$. We shall assume that $V$ decomposes into a direct sum of highest weight modules if restricted to 
the identity component $\widehat{{\cal L}\tilde{G}}^{k_f}$:
\begin{equation}
\label{decompositionreps}
V=\bigoplus_{\lambda\in I\subset\widehat{P}^+}V(\lambda )^{an}.
\end{equation}
Here, $I$ is a certain subset of those elements of $\widehat{P}^+$ having the same level.
Then , necessarily, the level $\lambda(c),\lambda\in I$ has to be a multiple of $k_f$.
Such representations are called negative energy representations of $\widehat{{\cal L}_{\Sigma}\tilde{G}}^{k_f}$.

Now, we want to construct all such representations of $\widehat{{\cal L}_{\Sigma}\tilde{G}}^{k_f}$.

Recall that $\Sigma \cong \widetilde{{\cal L}_{\Sigma}\tilde{G}}^{k_f}/\widetilde{{\cal L}\tilde{G}}^{k_f}$. 
We denote by $\widehat{P}^{+k}$ the set of dominant weights of level $k$.
Now, by a kind of Mackey induction we obtain the following result, see \cite{wendt:twistedcharacters}, Theorem 2.8 
and also \cite{toledanolaredo:representations}, Theorem 6.1 (Note that \cite{toledanolaredo:representations} considers
only projective representations):
\begin{theorem}
\label{nonconnirredreps}
Let $k$ be multiple of $k_f$. To each $\Sigma $-orbit $I\subset\widehat{P}^+$ there are $|\Sigma |/|I|$ many 
non-isomorphic irreducible representations of $\widetilde{{\cal L}_{\Sigma}\tilde{G}}^{k_f}$ whose restrictions
to $\widetilde{{\cal L}\tilde{G}}^{k_f}$ decompose according to Equation \ref{decompositionreps}.
\end{theorem}
We write $V(I)^{an}$ for the irreducible representation of $\widehat{{\cal L}\tilde{G}}^{k_f}$ given by $I\subset \widehat{P}^{+k}$
in this manner.
By the results of the previous subsection $\Sigma $  is generated by the element $\tau=\rho\sigma $ with $\rho\in Z$ and $\sigma\in\Gamma $.
Using Lemma \ref{descriptionauts} we can associate uniquely an element $\check{\beta}_{\rho }w_{\rho}\in \widehat{\W}(\check{P})_{\alc}$
representing the element $\rho\in\mbox{Aut}(\widehat{\Pi})$. As in the previous subsection we use the symbol $\check{\beta}_{\rho }$ 
also for the open loop $t\mapsto exp(2\pi i t \check{\beta}_{\rho })$. Choose a representative 
$n'_{w_{\rho}}\in N_{{\cal L}\tilde{G}}(T)$ of minimal possible order, where $T\subset\tilde{G}$ is a maximal torus of ${\cal L}\tilde{G}$.
Then it can easily be checked that the element $\check{\beta}_{\rho } n'_{w_{\rho}}\sigma $ is contained in 
${\cal L}_{\Sigma}\tilde{G}$. Denote by $n_{\tau}:=n_{\rho\sigma}$ one of 
its lifts to the central extension $\widetilde{{\cal L}_{\Sigma}\tilde{G}}^{k_f} $ preserving its order. 
This element clearly stabilises the maximal torus $\widetilde{T}$ of $\widetilde{{\cal L}\tilde{G}}^{k_f}$ corresponding to
$T$ and acts on $\widehat{T}:=\widetilde{T}\times\widetilde{\Cst}$ like the automorphism $\rho\sigma $. For Chevalley generators $e_{\widetilde{\alpha_i}},
f_{\widetilde{\alpha_i}},h_{\widetilde{\alpha_i}},\,\widetilde{\alpha_i}\in\widehat{\Pi} $ 
of $\widetilde{{\cal L}\liea}^{k_f}$ it is a straight forward calculation to check that $n_{\rho\sigma}$
permutes them according to its action on $\widehat{\Pi}$.
Furthermore, we find the following formula for the action of 
$n_{\rho\sigma}$ on $d$:
\begin{equation}
d\to d-\check{\beta}_{\rho }-\frac{k_f}{2}<\check{\beta}_{\rho },\check{\beta}_{\rho }>c.
\end{equation}

Now, consider the irreducible representation $V(I)^{an}$ where $n_{\rho\sigma}^{|I|}$ acts as the identity on 
$V(I)^{an}_{\lambda}$ for every highest weight vector $\lambda\in\widehat{P}^{+k}$. This determines the module $V(I)^{an}$ uniquely.
It is easy to see that we obtain a unitary action of $n_{\rho\sigma}$ on $V(\lambda)$.
Thus, there is a $\widehat{{\cal L}_{\Sigma}\tilde{G}}^{k_f}=\widetilde{{\cal L}_{\Sigma}\tilde{G}}^{k_f}\rtimes\widetilde{\Cst}$-action 
by requiring $d$ to act trivially on the highest weight spaces $V(I)^{an}$. 
(Recall that we have a $q$-fold cover $\widetilde{\Cst}\to\Cst $ which acts on $\widetilde{{\cal L}\tilde{G}}^{k_f} $.)
Denote its elements by $\tilde{q}$. 

As above we define the sets $\widehat{{\cal L}_{\Sigma}\tilde{G}}^{k_f}_{\tilde{q}}:=\widetilde{{\cal L}_{\Sigma}\tilde{G}}^{k_f}\times\{\tilde{q}\}$
and $\widehat{{\cal L}_{\Sigma}\tilde{G}}^{k_f}_{<1}:=\cup_{\tilde{q},\,|\tilde{q}|<1}\widehat{{\cal L}_{\Sigma}\tilde{G}}^{k_f}_{\tilde{q}}$.
Then we have the following result, see also \cite{wendt:twistedcharacters}, Corollary 2.9:
\begin{corollary}
Let $k$ be a multiple of $k_f$ and $I\subset\widehat{P}^{+k}$a $\Sigma$-orbit as above.
For any $\tilde{q}\in\tilde{\Cst}$ with $|\tilde{q}|<1$ and any $g\in\widetilde{{\cal L}_{\Sigma}\tilde{G}}^{k_f}$
the operator $g\tilde{q}^{-d}:V(I)^{an}\to V(I)^{an}$ uniquely extends to a trace class operator on 
$V(I)^{ss}:=\bigoplus_{\lambda\in I}V(\lambda)^{ss}$.
\end{corollary}
 
As before we define the character $\chi_{I}:\widehat{{\cal L}_{\Sigma}\tilde{G}}^{k_f}_{<1}\to\C$ for 
$(g,\tilde{q})\in \widehat{{\cal L}_{\Sigma}\tilde{G}}^{k_f}_{<1}$ by:
\begin{equation}
\chi_{I}(g,\tilde{q})=Tr_{V(I)^{ss}}(g\tilde{q}^{-d}).
\end{equation}
For $g\in \widetilde{{\cal L}\tilde{G}}^{k_f}=\widetilde{{\cal L}_{\Sigma}\tilde{G}}^{k_f}_0$ we clearly obtain 
$ \chi_{I}(g,\tilde{q})=\sum_{\lambda\in I} \chi_{\lambda}(g,\tilde{q})$.
Since $n_{\rho\sigma}$ acts as unitary operator the following result holds, see \cite{wendt:twistedcharacters}, Corollary 2.10:
\begin{corollary}
\label{convergentchars}
The functions $\chi_{I}$ are holomorphic and invariant under conjugation.
\end{corollary}

Our next aim is to give a more detailed description of the characters. We start by summarising the results on
conjugacy classes of the Kac-Moody group as outlined in \cite{wendt:twistedcharacters}, Section 3.
Let us fix the element $n_{\rho\sigma  }$ and the maximal torus
$\widehat{T}$ of the Kac-Moody group 
$\widehat{{\cal L}_{\Sigma}\tilde{G}}^{k_f}$ as above.
It is known that the Weyl group $\widehat{\cal W}$ can be described by the formula
$\widehat{\cal W}=N_{\widehat{{\cal L}_{\Sigma}\tilde{G}}^{k_f}_0}(\widehat{T})/\widehat{T}$. 
%respectively
%$\widehat{\cal W}=N_{\widehat{{\cal L}_{Z^{\sigma}}\tilde{G}}^{k_f}(\sigma)_0}(\widehat{T})/\widehat{T}$.
Consider the fixed point torus $\widehat{T}^{\rho\sigma}:=\{t\in\widehat{T}, \rho\sigma (t)=t\}$
which is easily seen to be connected. Let
$\widehat{{\cal L}_{\Sigma}\tilde{G}}^{k_f}_{\rho\sigma}$ be the connected component of the Kac-Moody group
%respectively $\widehat{{\cal L}_{Z^{\sigma}}\tilde{G}}^k(\sigma)_{\rho\sigma}$
corresponding to $\rho\sigma $. (Recall that $\rho\sigma $ generates the component group of our Kac-Moody group.)
Then, the shifted connected fixed point torus $\widehat{T}^{\rho\sigma}n_{\rho\sigma  }$ is contained in this connected component.
By results of Etingof etal., \cite{etingoffrenkelkirilov}, the set of all of its conjugates $\{gtg^{-1}, t\in \widehat{T}^{\rho\sigma}n_{\rho\sigma  },
g\in \widehat{{\cal L}_{\Sigma}\tilde{G}}^{k_f}_0\}$
is dense subset of this component.

Therefore, it is enough to determine the characters on $\widehat{T}^{\rho\sigma}_{<1}n_{\rho\sigma  }$.
Let us define the outer Weyl group $\widehat{\cal W}_{\rho\sigma  }:=
N_{\widehat{{\cal L}_{\Sigma}\tilde{G}}^{k_f}_0}(\widehat{T}^{\rho\sigma}_0n_{\rho\sigma})/\widehat{T}^{\rho\sigma}_0$ 
%respectively
%$\widehat{\cal W}=N_{\widehat{{\cal L}_{Z^{\sigma}}\tilde{G}}^{k_f}(\sigma)_0}(\widehat{T}^{\rho\sigma}_0\tau_{\rho\sigma})/\widehat{T}^{\rho\sigma}_0$
and the fixed point Weyl group $\widehat{\cal W}^{\rho\sigma  }$.

Then the following result holds, see \cite{wendt:twistedcharacters}, Section 3:
\begin{proposition}
(i) The outer Weyl group has the following structure:
\begin{equation}
\widehat{\cal W}_{\rho\sigma  }=\widehat{\cal W}^{\rho\sigma  }\ltimes(\widehat{T}/\widehat{T}^{\rho\sigma}_0)^{\rho\sigma}. 
\end{equation}
(ii) There is a natural $\widehat{\cal W}^{\rho\sigma  }$-equivariant isomorphism:
\begin{equation}  
\widehat{T}^{\rho\sigma}_0n_{\rho\sigma}/(\widehat{T}/\widehat{T}^{\rho\sigma}_0)^{\rho\sigma})\cong \widehat{T}/(1-\rho\sigma)(\widehat{T}).
\end{equation}
(iii) Two elements of $\widehat{T}^{\rho\sigma}_0\tau_{\rho\sigma}$ are conjugate under $\widehat{{\cal L}_{\Sigma}\tilde{G}}^{k_f}_0$ 
%respectively $\widehat{{\cal L}_{Z^{\sigma}}\tilde{G}}^{k_f}(\sigma)_0$ 
if and only if they are under $\widehat{\cal W}_{\rho\sigma  }$.
\end{proposition}
For the last result of this section we have to exclude the cases, where the folded Dynkin diagram is no longer of affine type, 
ie. the pairs $\mbox{A}_n^1,\rho$ with $\rho$ generating $\pi_1(PSl_{n+1})$.
Then, it is a simple observation that $\widehat{\cal W}^{\rho\sigma  }$ is the Weyl group of the orbit Lie algebra, the affine Kac Moody algebra 
corresponding to the folded Cartan matrix ${}^{\rho\sigma}\widehat{C}$. 
Furthermore, the co-invariant torus $\widehat{T}/(1-\rho\sigma)(\widehat{T})$ is 
$\widehat{\cal W}^{\rho\sigma  } $-equivariantly isomorphic to a maximal torus $\widehat{T}'$
the connected Kac-Moody group corresponding to the orbit Lie algebra.
Under this identification the lattice of $\Sigma$-invariant weights $\widehat{P}^{\Sigma}:=\{\lambda\in\widehat{P},\,\rho\sigma(\lambda)=\lambda\}$, 
corresponds to the 
weight lattice of the orbit Lie algebra, and this identification respects the cone of dominant weights.

By the second part of the proposition above, we can interpret the characters as functions on 
$(\widehat{T}/(1-\rho\sigma)(\widehat{T}))_{<1}:=\widehat{T}_{<1}/(1-\rho\sigma)(\widehat{T})$ which again is
$\widehat{T}'_{<1}$. Under the identifications above, 
the following theorem holds, see \cite{wendt:twistedcharacters}, Theorem 5.5 and \cite{fuchs}, Theorem 4.4.1: 
\begin{theorem}
Let $\lambda\in\widehat{P}^+$ a $\Sigma$-invariant dominant weight and denote by $\lambda'$ this same weight interpreted as an element 
of the orbit Lie algebra. Then, the corresponding characters coincide:
\begin{equation}
\chi_{\lambda}|_{\widehat{T}/(1-\rho\sigma)(\widehat{T}))_{<1}}=\chi'_{\lambda'}|_{\widehat{T}'_{<1}}.
\end{equation}
Here, $\chi'$ denotes the characters of the Kac-Moody group corresponding to the folded Dynkin diagram.
\end{theorem}
%%%
%%%
%%% Section 4: The $\C^{\ast}$-action on the section
%%%
%%%
\section{The Quotient Map and the Section}
%%%
%%%
%%%  Introductory part
%%%
%%%
This section contains our main result.
First, we will consider restriction of the conjugation action of 
the affine Kac Moody group $\widehat{{\cal L}_{\Sigma}\tilde{G}}^k$ respectively $\widehat{{\cal L}_{Z^{\sigma}}\tilde{G}}^k\!\!(\sigma)$
to an exterior component generating its component group. Using characters we shall define the `quotient map' on the $<1$-part in the
sense of Subsection 3.3. Then we prove our main result saying that this `quotient map' admits a section if we restrict it to 
$\tilde{q}\in\widetilde{\Cst}$ with $|\tilde{q}|$  
small. This section is defined as a transversal slice to a representative of twisted Coxeter element, (multiplied with $\tilde{q}$)
as introduced in section 1, and, furthermore, can be endowed with a $\Cst$-action. 

We keep the notation introduced in the previous sections. 
For the readability of the exposition the whole discussion will be carried out for the Kac Moody groups 
$\widehat{{\cal L}_{\Sigma}\tilde{G}}^{k_f}$.
However, all the results do hold for the twisted groups  $\widehat{{\cal L}_{Z^{\sigma}}\tilde{G}}^{k_f}\!\!(\sigma)$ as well.
As before, we write $\tau=\rho\sigma$ for the generator of its component group $\Sigma $.

We begin by defining the quotient map:
Consider the set of fundamental weights $\{\lambda_0,...,\lambda_r,\delta\}$. Then $\tau $ can be interpreted as a permutation of the index set 
$\{0,...,r\}$.
A set of dominant generators of the 
fixed point weight lattice $\widehat{P}^{\Sigma }$ is given by $\{\Lambda_0,...,\Lambda_{s-1},\delta\}$ where we use 
$\Lambda_i:=\frac{1}{|\Sigma_i|}\sum_{j=1}^{\tiny\mbox{ord}\,\tau } \lambda_{\tau^j(i)}$.

{\it Remark:} We have to be very careful because the $\Lambda_i$ do not coincide with the sum over $\Sigma $-orbit of 
$\lambda_i $. Indeed, a calculation yields $\tau (\lambda_i)-\lambda_{\tau(i)}\in\Q^{\ast}\delta $.

Taking a look at Table 5 one observes that the level of the irreducible $\kmliea $-module $V(\Lambda_i )^{an}$ is a multiple of the fundamental
level $k_f$ of ${\cal L}_{\Sigma}\tilde{G}$ and thus gives rise to an irreducible representation of 
$\widehat{{\cal L}_{\Sigma}\tilde{G}}^{k_f}$ on $V(\Lambda_i )^{an}$,
by Theorem \ref{nonconnirredreps}.
By Corollary \ref{convergentchars} the characters of these representations are holomorphic 
and conjugacy invariant functions on $\widehat{{\cal L}_{\Sigma}\tilde{G}}^{k_f}_{<1}$.

Motivated by the finite dimensional situation, see eg. \cite{steinberg:regularelements,ich:conjugacyclasses}, 
we define the `quotient map' (recall that $D^{\ast}$ is the punctured unit disk): 
\begin{eqnarray}
\chi:\widehat{{\cal L}_{\Sigma}\tilde{G}}^{k_f}_{<1, \tau}&\to & \C^s\times D^{\ast}\\
(g,\tilde{q})&\mapsto & (\chi_{\Lambda_0}(g,\tilde{q}),...,\chi_{\Lambda_{s-1}}(g,\tilde{q}),\chi_{\delta}(g,\tilde{q})=\tilde{q}).
\end{eqnarray}
Its restriction of $\chi $ to the set $\widehat{{\cal L}_{\Sigma}\tilde{G}}^{k_f}_{\tilde{q}, \tau }$ will be denoted by $\chi_{\tilde{q}} $.

We proceed by indicating the construction recipe for $S$:
Assume that $\Sigma $ acts on the set of simple roots $\widehat{\Pi} $ with $s$ orbits.
Choose a set $\{\alpha_0,...,\alpha_{s-1}\}$ of representatives of the $\Sigma $-orbits on $\widehat{\Pi}$.
For each index $i$ consider the corresponding canonical embedding $\phi_i:Sl_2\to  \widehat{{\cal L}_{\Sigma} G}^{k_f}$, the universal cover of $G$.
Let $X_{\alpha_i}:\C\to \widetilde{{\cal L}_{\Sigma}\tilde{G}}^{k_f}_0$ be the corresponding root group, ie 
$\mbox{im}X_{\alpha_i}=\phi_i(({\tiny \begin{array}{cc} 1 & \C \\ 0 & 1\end{array}}))$.
Furthermore, denote by $n_i\in N_{\widehat{{\cal L}_{\Sigma}\tilde{G}}^{k_f}_0}(\widehat{T})$ the following representative of 
the simple reflection $s_i\in\widehat{\W}$:
$n_i:= \phi_i(({\tiny \begin{array}{cc} 0 & -1 \\ 1 & 0\end{array}}))$.  
In addition, let $n_{\tau }\in N_{\widehat{{\cal L}_{\Sigma}\tilde{G}}^{k_f}}(\widehat{T})$ be the lift of $\tau\in\widehat{\cal W}\rtimes\Sigma$ 
defined in Subsection 3.3.
Now, we introduce the section $S$:
\begin{eqnarray}
S:\C^{s}\times \C^{\ast }&\to & \widehat{{\cal L}_{\Sigma}\tilde{G}}^{k_f}_{\tau,\,\tilde{q}} \\   
(c_0,...,c_{s-1},\tilde{q}) &\mapsto & (X_0(c_0)n_0...X_{s-1}(c_{s-1})n_{s-1}n_{\tau },\,\tilde{q}).
\end{eqnarray}
{\it Remark:} Note, that $S(0,...,0,1)$ is a representative in $N_{\widehat{{\cal L}_{\Sigma}\tilde{G}}^{k_f}}(\widehat{T})$
of the twisted Coxeter element.
Furthermore, we will write $S_{\tilde{q}}$ for the restriction $S|_{\C^{s}\times \{\tilde{q}\}}$. 

%Let us start by indicating the construction recipe for $S$:
%As above denote by $\Pi $ the set of simple roots of $\kmliea $ and let $\tau =\rho\sigma $ be an exterior automorphism of $\kmliea $.
%We will restrict to the case where $\tau $ is generator of the fundamental group $\pi_1(G)$ of a non simply-connected  
%algebraic group $G$ with Lie algebra $\liea $. 
%From now on we shall exclude the case where $G=PSl_n$, ie. $\tau $ is the symmetry of 
%$\mbox{A}^1_n$ shifting the cycle by one space.
%As before, denote by $\Sigma $ the subgroup of $\mbox{Aut }(G)$ generated by $\tau $.
%We will write $\widetilde{{\cal L} G}$ for a Kac-Moody group corresponding to $G$ with minimal possible level, cf. section 2.
%
Before stating our main result let us endow the image $\mbox{im}S_{\tilde{q}}$ with a $\Cst$-action:
We assume that the representatives $\{\alpha_0,...,\alpha_{s-1}\}$ are chosen in the same way as in the paragraph preceding Proposition \ref{descriptionb}.
Let $\iota:\Cst\to \widehat{{\cal L}_{\Sigma}\tilde{G}}^{k_f}_0$ be the parametrisation of the centre given by $c\in\kmliea$,
let $\mu_b\in\check{\chi}(\widehat{T})$ be the one parameter subgroup corresponding to $b\in\check{\widehat{Q}}$ 
appearing in Proposition \ref{descriptionb}, and let $k$ be the number given by $(\mbox{cox}^{\tau}-1)(b)=kc$.
We define the $\Cst$-action by:
\begin{eqnarray} 
\Cst\times \mbox{im}S_{\tilde{q}}& \to & \mbox{im} S_{\tilde{q}}\\
(u,S_{\tilde{q}}(c_0,...,c_{s-1}))&\mapsto & \iota (u)^{k} \mu_b (u)S_{\tilde{q}}(c_0,...,c_{s-1})\mu_b (u)^{-1}.
\end{eqnarray}
A simple calculation using Lemma \ref{twistedKaclabel} and Table 4 yields:
\begin{lemma}
\label{Cstarlemma}
For the $\Cst$-action we obtain explicitly:
\begin{equation}
\iota (u)^{k} \mu_b (u)S_{\tilde{q}}(c_0,...,c_{s-1})\mu_b (u)^{-1}=S_{\tilde{q}}(u^{k\Lambda_0(c)}c_0,...,u^{k\Lambda_{s-1}(c)}c_{s-1}).
\end{equation}
\end{lemma}
%%%
%%%
%%%   Theorem: Main result
%%%
%%%
Here comes our main result:
\begin{theorem}
\label{maintheorem}
For small $|\tilde{q}|$ the map $S_{\tilde{q}}$ is a section to $\chi_{\tilde{q}}$, ie. $\chi_{\tilde{q}}\circ S_{\tilde{q}}:\C^s\to\C^s$ is an isomorphism. 
\end{theorem}
The proof will basically proceed in the following way. First, using the explicit description of the cross section and
representation theory we will show that the Jacobian of the map $\chi\circ S$ is  a unit in the ring $\C\{\tilde{q}\}$,
the ring of convergent power series.
Then, the isomorphism follows by exploitation of the $\C^{\ast }$-action on the section and the quotient space. 

Before starting the proof of the theorem we compile some auxiliary results needed:
Let $V(\Lambda )$ be the irreducible $\kmliea _{\tiny\mbox{pol}}$-module of highest weight $\Lambda $
and assume $\Lambda$ to be $\Sigma$-invariant.
For the polynomial Lie algebra $\kmliea _{\tiny\mbox{pol}}$ we have a weight space decomposition:
\begin{equation}
V(\Lambda )=\bigoplus _{\mu \prec\Lambda}V(\Lambda )_{\mu },
\end{equation}
with weight spaces $V(\Lambda )_{\mu }$. As mentioned in subsection 3.3 the modules $V(\Lambda )^{an}$ are certain completions of these spaces, 
ie. we have to allow certain infinite series w.r.t. the above decomposition.

The following immediate lemma describes the action of certain group elements on the weight spaces:
%%%
%%%
%%%   Lemma 3.4: Action of group elements
%%%   actionlemma
%%%
\begin{lemma}
\label{actionlemma}
Let us keep the notation as above. Then, for all $v\in V(\lambda )_{\mu }^{an}$ and $t\in \widehat{T}$ one obtains:\\
i) $t.v=\mu (t) v$,\\
ii) $n_w.v\in V(\Lambda )_{w(\mu  )}^{an}$, for all $w\in \widehat{\W}$,\\
iii) $n_{\tau }.v\in V(\Lambda )^{an}_{\tau (\mu  )}$, if $\tau(\Lambda )=\Lambda$. \\
iv) $X_{\alpha }(c).v =v+\sum_{j=1}^{\infty }c^j v_j$ with $v_j\in V(\Lambda )^{an}_{\mu + j\alpha }$.
\end{lemma}
%%%
%%%
%%%  Lemma 3.5: Action of cross section.
%%%
%%%
For a weight $\mu $ of the representation $V(\Lambda )^{an}$ denote by $p_{\mu }$ resp. $i_{\mu }$ the canonical 
projections respectively injection of the corresponding weight space:
\begin{eqnarray}
i_{\mu }:V(\Lambda )^{an}_{\mu } & \hookrightarrow & V(\Lambda )^{an},\\
p_{\mu }:V(\Lambda )^{an} & \rightarrow & V(\Lambda )^{an}_{\mu }.
\end{eqnarray}

Let us investigate the action of elements of the cross section:
For calculating the characters evaluated on the section we are interested in those weights $\mu $ fulfilling:
\begin{equation}
p_{\mu }S(c_0,...,c_{s-1},\tilde{q})i_{\mu }\neq 0.
\end{equation}
Obviously , the following identity holds:
\begin{equation}
p_{\mu} S(c_0,...,c_{s-1},\tilde{q})i_{\mu} = \tilde{q}^{\mu (d)}p_{\mu} S(c_0,...,c_{s-1},1)i_{\mu}. 
\end{equation}  
There is a first result:
\begin{lemma}
For any weight $\mu $ of $V(\Lambda)^{an}$ with $\Lambda $ $\Sigma $-invariant the following holds:
\begin{equation}
p_{\mu }S(c_0,...,c_{s-1},\tilde{q})i_{\mu }=\tilde{q}^{\mu(d)}p_{\mu }X_{\alpha_0}(c_0)n_0 i_{\mu }...p_{\mu }X_{\alpha_{s-1}}(c_{s-1})n_{s-1}i_{\mu }p_{\mu } n_{\tau } i_{\mu }.
\end{equation}
Furthermore, this is only non-vanishing if $\mu\prec \Lambda $, $\tau (\mu )=\mu$ and $\mu $ dominant.
\end{lemma}
%%%
%%%
%%%   Proof of Lemma 3.5:
%%%   sectionlemma
%%%
{\it Proof:} First, note that
the following formula is true, see Br\"uchert \cite{bruechert}, Equation (35), and also \cite{ich:conjugacyclasses}, 
Equation (27), in the finite dimensional case:
\begin{equation}
\label{sectionlemma}
p_{\mu }X_{\alpha_0}(c_0)n_0...X_{\alpha_{s-1}}(c_{s-1})n_{s-1}i_{\mu }=
p_{\mu }X_{\alpha_0}(c_0)n_0 i_{\mu }...p_{\mu }X_{\alpha_{s-1}}(c_{s-1})n_{s-1}i_{\mu }. 
\end{equation}    
Now, the equation follows in the case $\tau (\mu )= \mu $.

The right hand side vanishes for $\tau (\mu )\neq\mu $ . For the left hand side we calculate for an element 
$v\in V(\Lambda )^{an}_{\mu }$ using Lemma \ref{actionlemma}:
\begin{equation}
X_{\alpha_0}(c_0)n_0...X_{\alpha_{s-1}}(c_{s-1})n_{s-1}n_{\tau }v \in 
\prod_{k_0,...,k_{s-1}\in\Z }V(\Lambda )^{an}_{\tau (\mu )+\sum_{i=0}^{s-1}k_i\alpha_i}.
\end{equation}
The equation and hence the second condition follow from the following observation:

{\sl Claim:} If $\tau (\mu )\neq \mu $ we have $\tau (\mu )-\mu \notin \Z<\alpha_0,...,\alpha_{s-1}>$.

Indeed, if we assume the contrary there exists $(k_0.,,,k_{s-1})\in\Z^s\backslash\{0\} $, such that 
$\tau (\mu )-\mu =\sum_{i=0}^{s-1}k_i\alpha_i$. Taking the sum over the $\Sigma $-orbit we would obtain:
$0=\sum_{j=0}^r l_j\alpha_j$ with $(l_0,...,l_r)\in\Z^r\backslash\{0\} $, which would be absurd. \hfill $\diamondsuit $.

Let us verify the non-vanishing conditions mentioned in the lemma. The first is automatic for highest weight representations
and we have already proven the second. 
For the third,  assume $\mu $ not to be dominant, ie 
$\mu=\sum_{j=0}^rn_j\Lambda_j +n\delta $ with $n_{i_0}<0$. Since $\tau(\mu )=\mu $ we may assume $0\leq i_0\leq s-1$.   
Then, $s_{i_0}\mu =\mu + n_{i_0}\alpha_{i_0}$ and because of Lemma \ref{actionlemma} (iv) we have 
$p_{\mu} S(c_0,...,c_{s-1},q)v=0$, for all $v\in V(\Lambda_j)_{\mu }$. \hfill $\qed $.

%%%
%%%
%%%  Corollary 3.1: Trace formula 
%%%  invariantcor
%%%
%\begin{corollary}
%\label{invariantcor}
%For the characters $\chi_{\Lambda_i}$ evaluated on  $tn_{\tau }\in \widehat{T}_{<1}n_{\tau }=\widehat{T}n_{\tau }\cap \widehat{{\cal L}_{\Sigma}\tilde{G}}^{k_f}_{<1}$ 
%we get:
%\begin{equation}
%\chi_{\Lambda_i}(tn_{\tau })=
%\sum_{\mu\prec\Lambda_i,\,\tau(\mu )=\mu }\mbox{dim}V(\Lambda_i )_{\mu } \mu (t) \mbox{Tr}\,(n_{\tau }|_{V(\Lambda_i )_{\mu }}).
%\end{equation}
%\end{corollary}
Our next goal is to determine those $\Sigma $-invariant weights $\mu $ satisfying the conditions mentioned in the Lemma.
%%%
%%%
%%%       Determination of contributing weights
%%%
%%%

By the following reasoning we get restrictions on the values of $n$ the coefficient of $\delta $
in $\mu =\sum_{j=0}^rn_j\Lambda_j +n\delta $:
Since the central variable $c$ has to act as a multiple of the identity on every irreducible representation,
we have $(\Lambda_i-\mu)(c)=0$.
This implies (for the stabiliser $\Sigma_i$ of $\alpha_i$ in $\Sigma$):
\begin{equation}
\sum_{l=1}^{\tiny\mbox{ord}\,\tau}\frac{\check{a}_{\tau^l(i)}}{|\Sigma_i|}=
\sum_{j=0}^{s-1}n_j\sum_{l=1}^{\tiny\mbox{ord}\,\tau}\frac{\check{a}_{\tau^l(j)}}{|\Sigma_i|},
\end{equation}
because $\delta (c)=0$ and  $\Lambda_i$ equals the sum over the 
$\Sigma $-orbit of certain $\Lambda_j$, up to a multiple of $\delta $. 
By the definition of folding of Dynkin diagrams we see that the sums 
$\sum_{l=1}^{\tiny\mbox{ord}\,\tau}\frac{\check{a}_{\tau^l(i)}}{|\Sigma_i|}$ coincide with the 
dual Coxeter labels ${}^{\tau}\check{a}_i$, $i\in\{0,...,s-1\}$, of the folded Cartan matrix ${}^{\tau}\widehat{C}$, for all $0\leq i\leq s-1$.
(Recall that we set ${}^{\tau }\check{a}_0=1$ in the case $\mbox{A}^1_{n},\gamma $ with $\gamma $ a generator of $\check{P}/\check{Q}$.)
%%%
%%%
%%% The sets D(\lambda)
%%%
%%%

Let us define the set ${\cal D}(\Lambda_i):=
\{\mu\in\widehat{P}^{\Sigma +} | V(\Lambda_i)_{\mu}^{an}\neq\{0\}\}$.
Calculating in $\widehat{P} /\Z\delta $ we have shown the statement below:
%%%
%%%
%%% Lemma 3.6: D(\lambda)mod\delta
%%% doflambdamoddeltalemma
%%%
\begin{lemma}
\label{doflambdamoddeltalemma}
Each element of ${\cal D}(\Lambda_i)\mbox{mod}\delta $ satisfies one of the following properties:\\
I: $\Lambda_i$\\
II: $\Lambda_j$ with $i\neq j$ and ${}^{\tau }\check{a}_i={}^{\tau }\check{a}_j$ or\\
III: $\sum_{j}n_j\Lambda_j$ with $\sum_{j}n_j\geq 2$ and $n_j=0$, if 
${}^{\tau }\check{a}_i\leq{}^{\tau }\check{a}_j$.\\
Furthermore, ${\cal D}(\Lambda_i\mbox{mod})\delta $ is finite.
\end{lemma}
%%%
%%%
%%% Definition of n_i(\mu)
%%%
%%%
%%%%%%%%%%%%%%%%%%%%%%% Warning: Check factor a_0!!!!!!!!!!!!
For $\mu\in {\cal D}(\Lambda_i)$ we set $n_i(\mu):=\mbox{max}\{n\in\Z|\,\mu+n\delta\in {\cal D}(\Lambda_i)$.
Then we can describe the set ${\cal D}(\Lambda_i)$ precisely:
%%%
%%%
%%%  Lemma 3.7: D(\lambda)
%%%  doflambdalemma
%%%
\begin{lemma}
\label{doflambdalemma}
Any element $\mu\in {\cal D}(\Lambda_i)$ has one of the following forms:\\
I: $\mu=\Lambda_i+n\delta$, $n\leq n_i(\Lambda_i)=0$\\
II: $\mu=\Lambda_j+n\delta$ with $i\neq j$ and ${}^{\tau }\check{a}_i={}^{\tau }\check{a}_j$, 
$n\leq n_i(\Lambda_i)\leq 0$  and\\
III: $\mu=\sum_{j}n_j\Lambda_j+n\delta$ with $\sum_{j}n_j\geq 2$ and $m_j=0$, if 
${}^{\tau }\check{a}_i\leq{}^{\tau }\check{a}_j$ and $n\leq n_i(\mu)\leq0$.
\end{lemma}
{\it Proof:}
By assumption there are non-negative integers $k_l,\,0\leq l\leq r$ satisfying:
\begin{equation}
\Lambda_i-\mu=\Lambda_i-(\sum_{j}n_j\Lambda_j+n\delta)=\sum_{l=0}^rk_l\alpha_l.
\end{equation}
Evaluation this formula on the derivation $d$ yields $-n=k_0$, whence the lemma.\hfill $\qed $.
 
%%%
%%%
%%%   Trace formula
%%%
%%%
Let us return to the characters evaluated on the section. 
Consider a $\Sigma $-invariant dominant integral weight $\mu=\sum_{i=0}^{s-1}m_i\Lambda_i+n\delta$. 
By Lemma \ref{actionlemma} there exits a linear map $\Phi^j_{\Lambda_i,\,\mu}\in\mbox{End}(V(\Lambda_i)_{\mu })$
with $p_{\mu }X_{\alpha_j}(c_j)n_j\iota_{\mu }=c_j^{m_j}\Phi^j_{\Lambda_i,\,\mu}$. 
We set $\Phi_{\Lambda_i,\,\mu}:=\Phi^0_{\Lambda_i,\,\mu}...\Phi^{s-1}_{\Lambda_i,\,\mu}n_{\tau}$.
This yields the following formula:
\begin{equation}
\chi_i(S(c_0,...,c_{s-1},\tilde{q}))=\sum_{\mu=\sum_{i=0}^{s-1}m_i\Lambda_i+n\delta,\,\mu\prec\Lambda_i} \tilde{q}^{-n}c_0^{m_0}...c_{s-1}^{m_{s-1}}
\mbox{Tr}|_{V(\Lambda_i)_{\mu}}\,\Phi_{\Lambda_i,\,\mu}.
\end{equation}
%%%
%%%
%%%  Analysis of trace formula
%%%
%%%
Applying Lemma \ref{doflambdalemma} yields $\chi_j\circ S\in\C[c_0,...,c_{s-1}]\{\tilde{q}\}$. Thus, for appropriately chosen $P\in\C[c_0,...,c_{s-1}]\{\tilde{q}\}$
we obtain:
\begin{equation}
\label{Eqchic}
\chi_i(S(c_0,...,c_{s-1},q))=a_i(\tilde{q}^{-1})c_i+\sum_{j\neq i,\,{}^{\tau }\check{a}_j={}^{\tau }\check{a}_i}a_{j\,i}(\tilde{q}^{-1})c_j+P(c_0,...,c_{s-1},\tilde{q})
\end{equation}
Lemma \ref{doflambdalemma} implies additionally that $P$ has degree at least two in the $c_0,...,c_{s-1}$ and does not depend on those $c_j $ with 
${}^{\tau }\check{a}_j\leq {}^{\tau }\check{a}_i$.

%%%
%%%
%%%  The Jacobian of \chi\circ S
%%%
%%%
Consider the Jacobian of the the map $\chi\circ S:\C^s\times\C^{\ast }\to \C^s\times\C^{\ast }$:
\begin{eqnarray}
{\cal J}\chi\circ S=
\left(
\begin{array}{cccc}
\frac{\partial \chi_0\circ S}{\partial c_0} & ... &\frac{\partial \chi_0\circ S}{\partial c_{s-1}} &
\frac{\partial \chi_0\circ S}{\partial \tilde{q}} \\
. & . & . & .\\ . & . & . & .\\ . & . & . & .\\
\frac{\partial \chi_{s-1}\circ S}{\partial c_0} & ... &\frac{\partial \chi_{s-1}\circ S}{\partial c_{s-1}} &
\frac{\partial \chi_{s-1}\circ S}{\partial \tilde{q}} \\
\frac{\partial \chi_{\delta}\circ S}{\partial c_0} & ... &\frac{\partial \chi_{\delta}\circ S}{\partial c_{s-1}} &
\frac{\partial \chi_{\delta}\circ S}{\partial \tilde{q}}
\end{array}
\right).
\end{eqnarray}
%%%
%%%
%%%  Its determinant
%%%
%%%
We are interested in its determinant. Since $ \chi_{\delta}\circ S(c_0,...,c_{s-1},\tilde{q})=\tilde{q}$ we get:
\begin{eqnarray}
\mbox{det}{\cal J}\chi\circ S=\mbox{det}
\left(
\begin{array}{ccc}
\frac{\partial \chi_0\circ S}{\partial c_0} & ... &\frac{\partial \chi_0\circ S}{\partial c_{s-1}}  \\
. & . & . \\ . & . & . \\ . & . & . \\
\frac{\partial \chi_{s-1}\circ S}{\partial c_0} & ... &\frac{\partial \chi_{s-1}\circ S}{\partial c_{s-1}} 
\end{array}
\right).
\end{eqnarray}
By reordering the index set $\{0,...,s-1\}$ in increasing order of the ${}^{\tau }\check{a_j}$ the matrix 
$(\frac{\partial \chi_i\circ S}{\partial c_j})$ will have upper block triangular shape:
\begin{eqnarray}
\left(\frac{\partial \chi_i\circ S}{\partial c_j}\right)=
\left(
\begin{array}{ccccc}
A_1 &&&&  \\
& A_2 & & \ast & \\ &  & .&& \\ & 0 & &.& \\
&&&& A_p
\end{array}
\right).
\end{eqnarray}
Thus we calculate for the determinant: $\mbox{det}{\cal J}\chi\circ S = \prod_{i=1}^p\mbox{det}A_i $.
By the shape of $\chi\circ S$ as described in Equation \ref{Eqchic} it only depends on $\tilde{q}$: 
$\mbox{det}{\cal J}\chi\circ S \in\C\{q\}$.

Concerning the blocks the following statement is true:
%%%
%%%
%%%  Proposition 3.4: determinant is unit.
%%%  detunitprop
%%%
\begin{proposition}
\label{detunitprop}
For all $i$ the determinant $\mbox{det}A_i$ is a unit in the ring $\C\{\tilde{q}\}$.
\end{proposition}

Before proving the proposition we need an auxiliary result:
By definition of the $X_{\alpha_i }$ and the $n_i $ we clearly have:
\begin{eqnarray}
X_{\alpha_j }(c_j)n_j|_{V(\Lambda_i)_{\Lambda_i}}=\mbox{id}_{V(\Lambda_i)_{\Lambda_i}}.
\end{eqnarray} 
Hence, $\Phi_{\Lambda_i,\,\Lambda_i}=p_{\Lambda_i}X_{\alpha_i}(c_i)n_i i_{\Lambda_i}$.
Then, the following lemma holds:
%%%
%%%
%%%  Lemma 3.8: Phi is not zero
%%%  diagonalunitlemma
%%%
\begin{lemma}
\label{diagonalunitlemma}
Let us keep the notation as above. Then:
\begin{equation}
\frac{\partial \chi_i\circ S}{\partial c_i}(c_0,...,c_{s-1},\tilde{q})=a_i(\tilde{q}),\quad\mbox{and}\quad a_i(0)\neq 0.
\end{equation}
\end{lemma}
{\it Proof:}
By definition of the $X_{\alpha_i }$ and the $n_i $ we clearly have:
\begin{eqnarray}
X_{\alpha_j }(c_j)n_j|_{V(\Lambda_i)_{\Lambda_i}}=\mbox{id}_{V(\Lambda_i)_{\Lambda_i}}.
\end{eqnarray} 
Hence, $\Phi_{\Lambda_i,\,\Lambda_i}=p_{\Lambda_i}X_{\alpha_i}(c_i)n_i i_{\Lambda_i}$. 
Thus, we only have to prove that $\Phi_{\Lambda_i,\,\Lambda_i}\neq 0$. Assume this not to be true.
Then, $X_{\alpha_i}(c_i)(V(\Lambda_i)_{\Lambda_i-\alpha_i})\subset V(\Lambda_i)_{\Lambda_i-\alpha_i}$.
Since $\Lambda_i+k\alpha_i$, $k>0$  is not a weight, $\Lambda_i-\alpha_i-k\alpha_i$ is not a weight, either.
Thus,  we must have $X_{-\alpha_i}(c_i)(V(\tilde{\lambda}_i)_{\tilde{\lambda}_i-\alpha_i})
\subset V(\Lambda_i)_{\Lambda_i-\alpha_i}$. This implies $n_i(V(\Lambda_i)_{\Lambda_i-\alpha_i})
\subset V(\Lambda_i)_{\Lambda_i-\alpha_i}$, because $n_i $ is contained in the subgroup generated by
$X_{\alpha_i}$ and $X_{-\alpha_i}$. this is a contradiction. \hfill $\qed $.

{\it Proof of the Proposition:} First, we prove the following claim:

{\sl Claim:}
Let $\Lambda_{j_1},...,\Lambda_{j_l}$, $l>1$ be distinct fundamental roots of the invariant lattice
$\widehat{P}^{\Sigma +}$ having the same dual Coxeter numbers, then $n_{j_i}(\Lambda_{j_{i+1}})\leq 0$ (with $l+1=1$) and 
this inequality is strict for at least one $i\in\{1,...,l\}$.

{\it Proof of Claim:} 
By Lemma \ref{doflambdalemma} we know that the inequality holds. To prove the strictness result we assume the contrary, ie.
$n_{j_i}(\Lambda_{j_{i+1}})= 0$, for all possible choices of the indices $j_1,..,j_l$ and all possible $l$.
By assumption there has to be sums of positive roots $\nu_i$ such that $\nu_i=\Lambda_{j_{i}}-
\Lambda_{j_{i+1}}\prec 0$.Taking the sum over all $i$ yields $\sum_i\nu_i=0$ implying $\nu_i=0$, for
all $i$ which is absurd. \hfill $\diamondsuit $.

Now, we proceed with the proof of the proposition: Using the definition of the determinant we get:
Assume $A_i$ to be a square matrix of size $l$.
\begin{equation}
\mbox{det}A_i=(A_i)_{11}...(A_i)_{ll}+\sum_{\sigma\in {\cal S}_l\backslash\{\mbox{id}\}}\mbox{sgn}\sigma 
\prod_{j=1}^l(A_i)_{j\sigma (j)}.
\end{equation}
By Lemma \ref{diagonalunitlemma} the first term on the right hand side is a unit in $\C\{\tilde{q}\}$. By the claim above and 
Lemma \ref{doflambdalemma} the remains sum on the right hand side is contained in $\tilde{q}\C\{\tilde{q}\}$.\hfill $\qed $.

Now we can proceed the proof of the main result:

{\it Proof of Theorem 3.1.:} By Proposition \ref{detunitprop} we know that $S_{\tilde{q}}$ is a local isomorphism for small $|\tilde{q}|$.
Now, observe that we $\chi_{\tilde{q}}\circ S_{\tilde{q}} :\C^s\to \C^s$ is equivariant with respect to the $\C^{\ast }$ actions 
on the image and the preimage space. On both spaces this action has the same weights with multiplicities which are all positive.
The result now follows from \cite{slodowy:singularities}, Section 8.1, Lemma 1,. \hfill $\qed $.
%%%
%%%
%%%  Appendix: Compilation of Tables
%%%
%%%
\appendix{\bf\Large Appendix}

In this appendix we compile the tables referred to in the text.
\begin{center}
\begin{tabular}{|c|c|c|}
type & $\tau $ & folded type \\
\hline 
$\mbox{A}^1_1$ & $\gamma $ &  $0$ \\
\hline
$\mbox{A}^1_n$ & $\gamma $ & $0$ \\
$n\, \mbox{even}$ & $\gamma^l,l|n+1,l>1 $ & $\mbox{A}^1_{l-1}$\\
 & $\sigma $ & $\mbox{A}^2_{n}$ \\
\hline
$\mbox{A}^1_n$ & $\gamma $ & $0$ \\
$n\, \mbox{odd}$ & $\gamma^l,l|n+1,l>1 $ & $\mbox{A}^1_{l-1}$\\
& $\sigma $ &  $\mbox{D}^2_{\frac{n+3}{2}}$\\
 & $\sigma\gamma $ &  $\mbox{C}^1_{\frac{n-1}{2}}$  \\
\hline
$\mbox{B}^1_n$ & $\gamma $ & $\mbox{A}^2_{2n-2}$   \\
\hline
$\mbox{C}^1_n, n\, \mbox{even}$ & $\gamma $ & $\mbox{A}^2_{2n}$   \\
\hline
$\mbox{C}^1_n, n\,\mbox{odd}$ & $\gamma $ & $\mbox{C}^1_{\frac{n-1}{2}}$  \\
\hline
$\mbox{D}^1_4$ & $\gamma $ & $\mbox{A}^2_2$ \\
 & $\gamma^2 $ & $\mbox{C}^1_2$  \\
 & $\sigma $ & $\mbox{A}^2_5$  \\
 & $\rho $ & $\mbox{D}^3_4$\\
\hline
$\mbox{D}^1_n$ & $\gamma $ & $\mbox{A}^2_{n-2}$  \\
$ n\, \mbox{even}$ & $\gamma^2 $ & $\mbox{C}^1_{n-2}$ \\
& $\sigma $ & $\mbox{A}^2_{2n-3} $    \\
& $\gamma\sigma $ & $\mbox{B}^1_{\frac{n}{2}}$  \\
\hline
$\mbox{D}^1_n$ & $\gamma $ & $\mbox{C}^1_{\frac{n-3}{2}}$  \\
$n\, \mbox{odd}$ & $\gamma^2 $ & $\mbox{C}^1_{n-2}$ \\
& $\sigma $ & $\mbox{A}^2_{2n-3} $    \\
& $\gamma\sigma $ & $\mbox{A}^2_{n-2}$  \\
\hline
$\mbox{E}^1_6$ & $\gamma $ & $\mbox{G}^1_2$  \\
& $\sigma $ &  $\mbox{E}^2_6$ \\
\hline
$\mbox{E}^1_7$ & $\gamma $ & $\mbox{F}^1_4$  \\
\hline
$\mbox{A}^2_{2n-1}$ & $\gamma $ & $\mbox{C}^1_{n-1}$   \\
\hline
$\mbox{D}^2_3$ & $\gamma $ &  $\mbox{A}^1_1$ \\
\hline
$\mbox{D}^2_{n+1}, n\, \mbox{even}$ & $\gamma $ & $\mbox{D}^2_{\frac{n}{2}+1}$  \\
\hline
$\mbox{D}^2_{n+1}, n\, \mbox{odd}$ & $\gamma $ & $\mbox{A}^2_{n-1}$\\ 
\hline
\end{tabular}
\end{center}
%\vspace{1ex}
\begin{center}
{\large T}ABLE 1: Folded Dynkin diagrams
\end{center}
\begin{center}
\begin{tabular}{|c|c|c|}
type & $\tau $ & $\chi(t)$ \\
\hline
$\mbox{A}^1_1$  & $\gamma $ & $(t-1)(t-1)$\\
\hline
$\mbox{A}^1_n,$ & $\gamma^l,\,l|n+1$ & $(t-1)(t^n-1)$\\
$n\,\mbox{even}$ & $\sigma$ & $(t-1)(t^n-1)$ \\
\hline
$\mbox{A}^1_n,$ & $\gamma^l,\,l|n+1$ & $(t-1)(t^n-1)$\\
$n\,\mbox{odd}$ & $\sigma$ & $(t^2-1)(t^{n-1}-1)$ \\
 & $\sigma\gamma $& $(t-1)(t^{n+1}-1)/(t+1)$\\
\hline
$\mbox{B}^1_n$ & $\gamma  $& $(t-1)(t^n-1)$\\
\hline
$\mbox{C}^1_n,\,n\,\mbox{even}$ & $\gamma $ & $(t^2-1)(t^{n-1}-1)$\\
\hline
$\mbox{C}^1_n,\,n\,\mbox{odd}$ & $\gamma $ & $(t-1)(t^{n}-1)$\\
\hline
$\mbox{D}^1_4$ & $\gamma $ & $(t^2-1)(t^3-1)$\\
 & $\gamma^2 $ & $(t-1)(t^4-1)$\\
 & $\sigma $& $(t^2-1)(t^3-1)$\\
 & $\rho $& $(t^2-1)(t^3-1)$\\
\hline
$\mbox{D}^1_n$ & $\gamma $ & $(t^2-1)(t^{n-1}-1) $ \\
$n\,\mbox{even}$ & $\gamma^2$ & $(t-1)(t^n-1)$\\
 & $\sigma $& $(t^2-1)(t^{n-1}-1)$\\
 & $\gamma\sigma $& $(t^4-1)(t^{n-3}-1)$\\
\hline
$\mbox{D}^1_n$ & $\gamma $ & $(t-1)(t^{n}-1) $ \\
$n\,\mbox{odd}$ & $\gamma^2$ & $(t-1)(t^n-1)$\\
 & $\sigma $& $(t^2-1)(t^{n-1}-1)$\\
 & $\gamma\sigma $& $(t-1)(t^2+1)(t^{n-2}-1)$\\
\hline
$\mbox{E}^1_6$ & $\gamma $ & $(t^2-1)(t^5-1)$\\
 & $\sigma $& $(t^3-1)(t^4-1)$\\
\hline 
$\mbox{E}^1_7$ & $\gamma $ & $(t^3-1)(t^5-1)$\\
\hline
$\mbox{A}^2_{2n-1}$ & $\gamma $ & $(t-1)(t^n-1)$\\
\hline
$\mbox{D}^2_{n+1}$ & $\gamma $ & $(t^2-1)(t^{n-1}-1)$\\ 
\hline
\end{tabular}
\end{center}
\vspace{1ex}
\begin{center}
{\large T}ABLE 2: Characteristic polynomials of twisted Coxeter elements
\end{center}
%%%
%%%
%%%  Table: values for b
%%%
%%%
\begin{center}
\begin{tabular}{|c|c|c|l|}
type &  $\tau $ & $k$  & $b$ \\
\hline 
$\mbox{A}^1_1$ & $\gamma $ & $1$ 
%& $(1,1)$ 
& $(1,0)$ \\
\hline
$\mbox{A}^1_n$ & $\gamma^l,l|n+1 $ & $l$ 
%& $(1,...,1)$ 
& $(n,n-1,...,1,0)$ \\
$ n\, \mbox{even}$ & $\sigma $ & $1$ 
%& $(1,...,1)$ 
& $(\frac{n}{2},\frac{n}{2},\frac{n}{2}-1,...,1,0,1,...,\frac{n}{2}-1)$ \\
\hline
$\mbox{A}^1_n$ & $\gamma^l,l|n+1 $ & $l$ 
%& $(1,...,1)$ 
& $(n,n-1,...,1,0)$ \\
$ n\, \mbox{odd}$ & $\sigma $ & $1$ 
%& $(1,...,1)$ 
&  $(\frac{n-1}{2},\frac{n-1}{2},\frac{n-3}{2},...,1,0,1,...,\frac{n-3}{2})$ \\
$$ & $\sigma\gamma $ & $1$
%& $1$ & $(1,...,1)$ 
&  $(\frac{n+1}{2},\frac{n-1}{2},\frac{n-3}{2},...,1,0,1,...,\frac{n-1}{2})$ \\
\hline
$\mbox{B}^1_n, n\, \mbox{even}$ & $\gamma $ & $1$ 
%& $(1,1,2,...,2,1)$ 
&  $(\frac{n}{2}-1,\frac{n}{2},n-2,n-3,...,1,0)$ \\
\hline
$\mbox{B}^1_n, n\, \mbox{odd}$ & $\gamma $ & $2$ 
%& $(1,1,2,...,2,1)$ 
&  $(n-2,n,2n-4,2n-6,...,2,0)$ \\
\hline
$\mbox{C}^1_n, n\, \mbox{even}$ & $\gamma $ & $1$ 
%& $(1,...,1)$ 
&  $(\frac{n}{2},\frac{n}{2}-1,\frac{n}{2}-2,...,1,0,0,1,...,\frac{n}{2}-1)$ \\
\hline
$\mbox{C}^1_n, n\,\mbox{odd}$ & $\gamma $ & $1$ 
%& $(1,...,1)$ 
&  $(\frac{n+1}{2},\frac{n-1}{2},\frac{n-3}{2},...,1,0,1,...,\frac{n-1}{2})$ \\
\hline
$\mbox{D}^1_4$ & $\gamma $ & $1$ 
%& $(1,1,2,1,1)$ 
& $(3,1,2,0,2)$ \\
$$ & $\gamma^2 $ & $1$ 
%& $(1,1,2,1,1)$ 
& $(1,2,2,1,0)$ \\
$$ & $\sigma $ & $2$ 
%& $(1,1,2,1,1)$ 
& $(3,3,4,2,0)$ \\
$$ & $\rho $ & $1$ 
%& $(1,1,2,1,1)$ 
& $(1,2,1,1,0)$ \\
\hline
$\mbox{D}^1_n$ & $\gamma $ & $1$ 
%& $(1,1,2,...,2,1,1)$ 
& $(\frac{n}{2}-2,\frac{n}{2},n-4,...,2,0,2,...,n-6,\frac{n}{2}-1,\frac{n}{2}-3)$ \\
$ n\, \mbox{even}$ & $\gamma^2 $ & $1$ 
%& $(1,1,2,...,2,1,1)$ 
& $(\frac{n}{2}-1,\frac{n}{2},n-2,...,2,1,0)$ \\
$$ & $\sigma $ & $2$ 
%& $(1,1,2,...,2,1,1)$ 
& $(n-1,n-1,2n-4,2n-6,...,4,2,0)$ \\
 & $\gamma\sigma $ & $1$ 
%& $(1,1,2,...,2,1,1)$ 
& $(\frac{n}{2}-1,\frac{n}{2}-1,n-4,...2,0,2,...,\frac{n}{2}-2,\frac{n}{2}-2)$ \\
\hline
$\mbox{D}^1_n$ & $\gamma $ & $1$ 
%& $(1,1,2,...,2,1,1)$ 
& $(\frac{n-3}{2},\frac{n+1}{2},n-3,...,2,0,2,...,n-5,\frac{n-1}{2},\frac{n-5}{2})$ \\
$ n\, \mbox{odd}$ & $\gamma^2 $ & $2$ 
%& $(1,1,2,...,2,1,1)$ 
& $(n-2,n,2n-4,...,4,2,0)$ \\
$$ & $\sigma $ & $1$ 
%& $(1,1,2,...,2,1,1)$ 
& $(\frac{n-1}{2},\frac{n-1}{2},n-2,n,n-3,...,2,1,0)$ \\
$\mbox{D}^1_n, n\, \mbox{odd}$ & $\gamma\sigma $ & $1$ 
%& $(1,1,2,...,2,1,1)$ 
& $(\frac{n-1}{2},\frac{n-1}{2},n-3,...2,0,2,...,\frac{n-3}{2},\frac{n-3}{2})$ \\
\hline
$\mbox{E}^1_6$ & $\gamma $ & $2$ 
%& $(1,1,2,3,2,1,2)$ 
& $(3,7,8,3,4,5,0)$ \\
 & $\sigma $ & $2$ 
%& $(1,1,2,3,2,1,2)$ 
& $(1,5,6,3,,2,3,0)$ \\
\hline
$\mbox{E}^1_7$ & $\gamma $ & $2$ 
%& $(1,2,3,4,3,2,1,2)$ 
& $(7,10,9,4,3,6,5,0)$ \\
\hline
$\mbox{A}^2_{2n-1}, n\, \mbox{even}$ & $\gamma $ & $1$ 
%& $(1,1,2,...,2)$ 
&  $(\frac{n}{2}-1,\frac{n}{2},n-2,...,1,0)$ \\
\hline
$\mbox{A}^2_{2n-1}, n\, \mbox{odd}$ & $\gamma $ & $2$ 
%& $(1,1,2...,2)$ 
&  $(n-2,n,2n-4,...,2,0)$ \\
\hline
$\mbox{D}^2_{n+1}, n\, \mbox{even}$ & $\gamma $ & $1$ 
%& $(1,2,...,2,1)$ 
&  $(\frac{n}{2},n-2,...,2,0,0,2,...,n-4,\frac{n}{2}-1)$ \\
\hline
$\mbox{D}^2_{n+1}, n\, \mbox{odd}$ & $\gamma $ & $1$ 
%& $(1,2,...,2,1)$ 
&  $(\frac{n+1}{2},n-1,...,2,0,2,...,n-3,\frac{n-1}{2})$\\
\hline 
\end{tabular}
\end{center}
\vspace{1ex}
\begin{center}
{\large T}ABLE 3: Values for $b$ from Proposition 2.9
\end{center}
\begin{center}
\begin{tabular}{|c|c|c|}
type & $\tau$ & $p$  \\
\hline
$\mbox{A}^1_1$ & $\gamma $ & $2$\\
\hline
$\mbox{A}^1_n$ & $\gamma^l,l|n+1 $ & $n+1$\\
$ n\, \mbox{even}$ & $\sigma $ & $1$ \\
\hline
$\mbox{A}^1_n$ & $\gamma^l,l|n+1 $ & $n+1$\\
$ n\, \mbox{odd}$ & $\sigma $ &  $1$\\
 & $\sigma\gamma $ &  $2$  \\
\hline
$\mbox{B}^1_n, n\, \mbox{even}$ & $\gamma $ & $1$   \\
\hline
$\mbox{B}^1_n, n\, \mbox{odd}$ & $\gamma $ & $2$   \\
\hline
$\mbox{C}^1_n, n\, \mbox{even}$ & $\gamma $ & $1$   \\
\hline
$\mbox{C}^1_n, n\,\mbox{odd}$ & $\gamma $ & $2$  \\
\hline
$\mbox{D}^1_4$ & $\gamma $ & $2$ \\
 & $\gamma^2 $ & $2$  \\
 & $\sigma $ & $2$  \\
 & $\rho $ & $1$\\
\hline
$\mbox{D}^1_n$ &  $\gamma $ & $2$  \\
$ n\, \mbox{even}$ & $\gamma^2 $ & $4$ \\
 & $\sigma $ & $2$    \\
 & $\gamma\sigma $ & $2$  \\
\hline
$\mbox{D}^1_n$ & $\gamma $ & $4$  \\
$n\, \mbox{odd}$ & $\gamma^2 $ & $4$ \\
 & $\sigma $ & $1$    \\
 & $\gamma\sigma $ & $2$  \\
\hline
$\mbox{E}^1_6$ & $\gamma $ & $6$  \\
 & $\sigma $ &  $2$ \\
\hline
$\mbox{E}^1_7$ & $\gamma $ & $4$  \\
\hline
$\mbox{A}^2_{2n-1}, n\, \mbox{even}$ & $\gamma $ & $2$   \\
\hline
$\mbox{A}^2_{2n-1}, n\, \mbox{odd}$ & $\gamma $ & $4$   \\
\hline
$\mbox{D}^2_{n+1}$ & $\gamma $ & $2$  \\
\hline
\end{tabular}
\end{center}
\vspace{1ex}
\begin{center}
{\large T}ABLE 4: The factor $p$ of Lemma 2.10
\end{center}
\begin{center}
\begin{tabular}{|c|c|c|}
type & $\tau$ & $k_f$  \\
\hline
$\mbox{A}^1_1$ & $\gamma $ & $2$\\
\hline
$\mbox{A}^1_n$ & $\gamma^l,l|n+1 $ & $1$\\
$n\, \mbox{even}$& $\sigma $ & $1$ \\
\hline
$\mbox{A}^1_n$& $\gamma^l,l|n+1,\,l\,\mbox{even} $ & $1$\\
$n\, \mbox{odd}$ & $\gamma^l,l|n+1,\,l\,\mbox{odd} $ & $2$\\
 & $\sigma $ &  $1$ \\
 & $\sigma\gamma $ &  $1$ \\
\hline
$\mbox{B}^1_n$ & $\gamma $ & $1$  \\
\hline
$\mbox{C}^1_n, n\, \mbox{even}$ & $\gamma $ & $1$  \\
\hline
$\mbox{C}^1_n, n\,\mbox{odd}$ & $\gamma $ & $2$  \\
\hline
$\mbox{D}^1_4$ & $\gamma $ & $2$  \\
 & $\gamma^2 $ & $1$  \\
 & $\sigma $ & $1$ \\
 & $\rho $ & $1$  \\
\hline
$\mbox{D}^1_n$ & $\gamma $ & $2$  \\
$ n\, \mbox{even}$& $\gamma^2 $ & $1$ \\
& $\sigma $ & $1$   \\
$n\equiv 0 (4)$& $\gamma\sigma $ & $1$   \\
$ n\equiv 2 (4)$ & $\gamma\sigma $ & $2$   \\
\hline
$\mbox{D}^1_n$ & $\gamma $ & $2$  \\
$n\, \mbox{odd} $ & $\gamma^2 $ & $1$ \\
 & $\sigma $ & $1$   \\
 & $\gamma\sigma $ & $1$  \\
\hline
$\mbox{E}^1_6$ & $\gamma $ & $1$  \\
 & $\sigma $ &  $1$ \\
\hline
$\mbox{E}^1_7$ & $\gamma $ & $2$    \\
\hline
$\mbox{A}^2_{2n+1}$ & $\gamma $ & $1$   \\
\hline
$\mbox{D}^2_{n+1}$ & $\gamma $ & $1$ \\
\hline 
\end{tabular}
\end{center}
\vspace{1ex}
\begin{center}
{\large T}ABLE 5: Fundamental levels
\end{center}

%%%
%%%
%%% End of document, bibliography
%%%
%%%
%\bibliographystyle{plain}
%\bibliography{index}

\end{document}